\documentclass[preprint, 3p, 11pt, authoryear]{elsarticle}
\usepackage{graphicx}
\usepackage{amssymb}
\usepackage[utf8]{inputenc}
\usepackage{algorithm}
\usepackage{algorithmicx}
\usepackage{algpseudocode}
\usepackage{amsmath}
\usepackage{amsthm}
\usepackage{dsfont}
\usepackage{enumitem}
\usepackage{etaremune}
\usepackage{eurosym}
\usepackage{natbib}
\usepackage{listings}
\usepackage{xcolor}
\usepackage{bm}
\usepackage{booktabs}
\usepackage{nicefrac}
\usepackage{graphicx}
\usepackage{multirow}
\usepackage{color, colortbl}
\usepackage{xcolor}
\definecolor{Gray}{gray}{0.9}
\usepackage{url}
\usepackage{caption}
\usepackage{subcaption}
\usepackage{enumitem}
\usepackage{microtype} 
\usepackage{multicol}
\usepackage{wrapfig}
\usepackage{framed} 
\usepackage{nomencl}
\usepackage{tablefootnote}

\makenomenclature
\renewcommand*\nompreamble{\begin{multicols}{2}}
\renewcommand*\nompostamble{\end{multicols}}
\newtheoremstyle{break}
 {\topsep}{\topsep}%
 {\itshape}{}%
 {\bfseries}{}%
 {\newline}{}%
\theoremstyle{break}

\newcommand{\beginsupplement}{%
 \setcounter{table}{0}
 \renewcommand{\thetable}{S\arabic{table}}%
 \setcounter{figure}{0}
 \renewcommand{\thefigure}{S\arabic{figure}}%
 }

\journal{European Journal of Operations Research}

\begin{document}

\newtheorem{definition}{Definition}
\newtheorem{theorem}{Theorem}
\begin{frontmatter}


\title{Scalable Policies for \\ the Dynamic Traveling Multi-Maintainer Problem with Alerts}

% alphabetical order (surnames):

\author[sor,tue]{Peter Verleijsdonk\corref{corauthor}}
\ead{p.verleijsdonk@tue.nl}
\author[ieis,tue]{Willem van Jaarsveld}
\ead{w.l.v.jaarsveld@tue.nl} 
\author[sor,tue]{Stella Kapodistria}
\ead{s.kapodistria@tue.nl}

\address[sor]{Department of Mathematics and Computer Science}
\address[ieis]{Department of Industrial Engineering and Innovation Sciences}
\address[tue]{Eindhoven University of Technology, P.O. Box 513, Eindhoven 5600 MB, the Netherlands}
\date{\today}

\cortext[corauthor]{Corresponding author}

\begin{abstract}
Downtime of industrial assets such as wind turbines and medical imaging devices is costly. To avoid such downtime costs, companies seek to initiate maintenance just before failure, which is challenging because: (i) Asset failures are notoriously difficult to predict, even in the presence of real-time monitoring devices which signal degradation; and (ii) Limited resources are available to serve a network of geographically dispersed assets. In this work, we study the dynamic traveling multi-maintainer problem with alerts ($K$-DTMPA) under perfect condition information with the objective to devise scalable solution approaches to maintain large networks with $K$ maintenance engineers. Since such large-scale $K$-DTMPA instances are computationally intractable, we propose an iterative deep reinforcement learning (DRL) algorithm optimizing long-term discounted maintenance costs. The efficiency of the DRL approach is vastly improved by a reformulation of the action space (which relies on the Markov structure of the underlying problem) and by choosing a smart, suitable initial solution. The initial solution is created by extending existing heuristics with a dispatching mechanism. These extensions further serve as compelling benchmarks for tailored instances. We demonstrate through extensive numerical experiments that DRL can solve single maintainer instances up to optimality, regardless of the chosen initial solution. Experiments with hospital networks containing up to $35$ assets show that the proposed DRL algorithm is scalable. Lastly, the trained policies are shown to be robust against network modifications such as removing an asset or an engineer \underline{or yield a suitable initial solution} for the DRL approach.
\end{abstract}

\begin{keyword}
Maintenance \sep Dispatching \sep Decision process \sep Traveling maintainer problem \sep Deep reinforcement learning


\end{keyword}
\end{frontmatter}

% Load tex files
\section{Introduction}\label{sec:introduction}
Industrial assets such as medical imaging equipment, wind turbines and wafer steppers are expected to operate on demand. To ensure maximum availability, assets nowadays are regularly inspected to evaluate their degradation, i.e., \textit{condition-based maintenance} (CBM). Equipping assets with sensory equipment for real-time degradation monitoring can be used to devise efficient CBM policies. For example, an array of sensors is connected to a magnetic resonance imaging device that tracks medical procedures, but the gathered data can also be used to signal various stages of degradation and emit alerts, which in turn may be used to dispatch a maintenance engineer. Decision-makers often rely on heuristics or human intuition which often decouple dispatching from other operational decisions, for instance, the dispatching and relocation decisions of emergency service providers after an incident occurred. Devising maintenance and dispatching policies are, each in their own right, notoriously difficult problems and have been studied individually. 

The combination of maintenance and dispatching decisions, although understudied, is extremely important. As an essential step in the direction of the application, we contribute to the literature with the Dynamic Traveling Multi-Maintainer Problem with Alerts ($K$-DTMPA). In the $K$-DTMPA, $K$ engineers travel in a network of assets where the degradation of each asset is modeled through a finite number of degradation states. The first degradation state captures that the asset is as-good-as-new. Subsequently, the severity of degradation increases with the degradation state until the asset reaches the failed state and becomes unavailable. After each state transition, an alert is issued immediately to a central decision-maker who is responsible for scheduling maintenance and dispatching the maintenance engineers. We propose heuristic dispatching solutions based on classical ranking heuristics, but the optimal solution likely cannot be represented by a simple set of rules. Indeed, the $K$-DTMPA model can be naturally formulated as a Markov decision process (MDP), which is computationally intractable for realistic size problems. It is well known that such realistic sequential decision-making problems suffer from the curse of dimensionality.

This curse of dimensionality can be tackled using approximate dynamic programming / deep reinforcement learning (DRL) \citep{boute2021deep,powell2019unified} via a combination of machine learning and simulation. DRL has achieved state-of-the-art performance, for instance, in Atari video games \citep{atari} and chess \citep{chess}. However, while numerical experiments for sequential decision-making problems arising in operations management have yielded encouraging results, the ability to solve industrial-scale instances is restricted due to the long training times \citep{inventory_management}. Indeed, DRL has been shown to produce near-optimal policies for Dynamic Traveling Maintainer Problem with Alerts (DTMPA) instances \citep{dtmpa}, but training times exceed 12 hours even for networks containing 6 assets maintained by a single engineer. The $K$-DTMPA instances studied here are more challenging because they involve \emph{multiple} engineers that need to be dispatched in a coordinated fashion to up to 35 assets, and overcoming this challenge is a crucial step towards industrial-scale solutions. 

We adopt a form of approximate policy iteration (API) - an iterative algorithm in the DRL domain - and demonstrate that heuristic solutions may be leveraged as a starting point for the training algorithm, and that doing so may vastly reduce training times. More specifically, we extend the ranking heuristics for the DTMPA framework proposed by \cite{dtmpa} to the $K$-DTMPA framework by equipping them with a state-of-the-art dispatching heuristic. We demonstrate that API can improve such heuristics and can learn a repositioning strategy for unassigned maintenance engineers aimed at anticipating future alerts and failures. Additionally, to significantly reduce the action space complexity, we propose a suitable reformulation of the associated MDP in which actions for engineers are selected sequentially in a fixed order. Using an engineer-centric feature representation for this MDP reformulation further improves DRL's efficiency. 

The $K$-DTMPA represents a rich class of problems, and as a consequence, it is challenging to design general-purpose traditional heuristic algorithms that can be used to benchmark our learned policies: Such heuristics would have to incorporate jointly the geographical layout, the observed degradation, the costs for asset unavailability, maintenance and travel, and the spatial and temporal information regarding the engineers. To circumvent this, we devise two suitable subclasses of $K$-DTMPA instances that allow for the construction of strong benchmarks, namely the \emph{single maintainer} and the \emph{dispatching \& repositioning} (D\&R) instances. 
We show that API can train policies that outperform the benchmark for such instances. Moreover, our algorithm produces state-of-the-art policies for more complex instances.

The primary contributions of the paper are specified as follows:
\begin{itemize}
 \item The proposed $K$-DTMPA model jointly optimizes maintenance and dispatching decisions, problems that are linked in practice but are studied separately in prior literature. 
\item To reduce the action and state space complexity, we propose an MDP reformulation in which actions for the engineers are selected sequentially using an engineer-centric feature representation, and show that this yields more cost-effective $K$-DTMPA policies.
 \item We propose a generic approach that leverages classical heuristic policies to improve the training of neural network policies, and we demonstrate its effectiveness for $K$-DTMPA instances.
\end{itemize}

The main insights gathered from the numerical experiments are:
\begin{itemize}
 \item API can solve single maintainer instances up to optimality within a few iterations, regardless of the choice of the initial solution.
 \item Sophisticated dispatching heuristics are superior initial solutions when solving multi-maintainer instances compared to trivial policies such as the random policy.
 \item The trained policies are robust against removing an asset/engineer or yield suitable initial solutions to optimize such instances.
\end{itemize}

 By providing scalable solution approaches to make data-driven decisions for industrial-scale problems, we attempt to bridge the gap between academia and industry.

The remainder of the paper is structured as follows. \mbox{Section \ref{sec:literature}} provides an overview of the related literature. \mbox{Section \ref{sec:model}} formalizes in detail the $K$-DTMPA framework. In \mbox{Sections \ref{sec:heuristics} and \ref{sec:api}}, we detail the heuristic solutions and the deep reinforcement learning algorithm, respectively. In \mbox{Section \ref{sec:numerical_experiments}}, the setup of the numerical experiments on Dutch hospital networks is presented. \mbox{Section \ref{sec:numerical_results}} discusses the numerical results and the corresponding managerial insights. We conclude in \mbox{Section \ref{sec:conclusion}} and discuss operational constraints limiting application in industry.
\section{Literature review}\label{sec:literature}

In this section, we first discuss relevant literature in the streams of maintenance optimization and traveling maintainer problems. Subsequently, we provide an overview of the application of DRL in dynamic dispatching problems integrating operational decisions with a focus on scalable DRL for decision-making.

\subsection*{Maintenance optimization and traveling maintainer problems}\label{subsec:maintenance_opt_TMP}
According to \citet{OLDEKEIZER2017405,DEJONGE2020805}, maintenance models can be single-asset or multi-asset. Multi-asset models generalize single-asset models by considering joint maintenance policies for assets with any of the following dependencies: economic, structural, stochastic or resource dependency \citep{OLDEKEIZER2017405}. The degradation of assets is often modeled using a stochastic process that takes values in a discrete finite state-space, e.g., a Markov chain. Scheduled inspections can be improved by leveraging information acquired via sensors. These sensors sometimes measure asset degradation directly, for instance in the form of alerts \citep{DEJONGE201693, AKCAY2021}. \citet{abdul2018optimally} study the problem of optimally replacing multiple stochastically degrading systems using condition-based maintenance in a shared environment. However, in the multi-asset scenario, the geographical layout often constitutes a complex dependency, prompting the traveling maintainer problem (TMP). The goal of the traditional TMP is to find a route that visits each asset such that the sum of the times needed to reach each asset is minimized. The TMP is a mean-flow variant of the traveling salesman problem (TSP) and is thus NP-complete \citep{afrati1986complexity}. The computational complexity further increases when assigning a hard deadline to each asset, e.g., a bound on the response time. The TMP objective of minimizing the sum of functions of response times is studied by \citet{camci2014travelling}. Real-time CBM prognostics are incorporated in a TSP by scheduling maintenance using forecasted failure information, which is then generalized to also include travel times \citep{camci2015maintenance}. The dynamic TMP (DTMP) considers jobs that appear uniformly in a region according to a Poisson process \citep{bertsimas1989dynamic}. These jobs must be completed by a single maintainer with the objective of minimizing the average response time. \citet{drent2020dynamic} model a DTMP variant as a sequential decision-making problem and provide heuristic solution approaches leveraging real-time condition information to the dispatching and repositioning subproblems based on the minimum weighted bipartite matching problem and the maximum expected covering location problem, respectively. \citet{pechina2019real} propose and evaluate a range of heuristic solution approaches to the dispatching and relocation subproblems inspired by the domain of emergency response networks to serve a network of identical geographically distributed assets. In emergency response dispatching problems, it is commonly believed that the closest idle ambulance rule is near-optimal, however, significant cost reductions can be achieved by dispatching policies that account for coverage \citep{jagtenberg} and relocation of ambulances \citep{jagtenberg2}. 

Condition-based maintenance optimization typically optimizes the timing of maintenance, taking into account risks, costs and dependencies \citep{DEJONGE201693}, which is a formidable problem in itself. In such models, the dispatching of resources based on their spatio-temporal availability is rarely modeled in detail, and resources are typically abstracted away. To the best of our knowledge, dispatching and coordinating resources in response to unforeseen alerts has only been studied outside of the (condition-based) maintenance context, e.g., in ambulance dispatching \citep{jagtenberg2}. Our newly proposed $K$-DTMPA model jointly considers resource dispatching and tactical postponement of maintenance, which advances the applicability of condition-based maintenance models in a practical context, while in a sense merging two streams of literature.

\subsection*{The application of deep reinforcement learning in dynamic dispatching problems}\label{subsec:DRL_dynamic_dispatching}
Recent advances in machine learning have led to a variety of applications in various fields. For example, applications in the field of dynamic dispatching include ambulance dispatching, ATM servicing and mining logistics. Relevant challenges of the application of DRL include: (i) multi-agent systems may have a variable number of agents, (ii) variable objectives require costly retraining, (iii) the curse of dimensionality, (iv) the stochastic environment can be non-stationary and (v) the explainability of the trained agent \citep{khorasgani2020challenges}. \citet{zhang2020dynamic} tackle the Open-Pit Operational Planning problem by training a neural network that is shared amongst all the agents (trucks), i.e., the network receives each agent’s observation and outputs actions for each agent independently. \citet{holler2019deep} also apply DRL from a system-centric perspective to solve the multi-driver vehicle D\&R problem. \citet{schmid2012solving} solves the dynamic ambulance D\&R problem using approximate dynamic programming on a Vienna case study. \citet{ji2019deep} provide an effective dynamic ambulance redeployment algorithm implementing a neural network trained to score the waiting locations. \Citet{dtmpa} integrate maintenance and dispatching decisions in a holistic DTMPA framework, including uncertainty in the acquired information in the form of three information levels. They propose a wide range of heuristic solution approaches and a DRL algorithm to optimize long-term discounted costs.

DRL has the potential to deliver good policies for any operations management problem that possesses a natural MDP formulation. MDP instances considered in DRL studies are often restricted to stylized models, which is in contrast with the complexity of practical problems arising in operations management \citep{boute2021deep}. Approaches to make DRL more scalable include aggregating states \citep{StateAggregation} or modifying action selection, e.g., decoupling action selection \citep{ride_hailing} or by using continuous action representations \citep{Vanvuchelen2022}. 
\citet{DEMOOR2022535} and \citet{boute2021deep} argue that incorporating domain knowledge embedded in well-performing heuristic policies into the training algorithm improves DRL's efficiency, e.g., through \emph{reward shaping}. Reward shaping incentivizes the DRL agent to act similar to the action selected by the heuristic policy.

Like \citet{DEMOOR2022535}, we propose an approach that leverages domain knowledge to improve DRL's efficiency. Structurally, our approach differs from the reward shaping approach adopted by \citet{DEMOOR2022535} as follows: Reward shaping \emph{alters} the MDP formulation to reward actions that coincide with the actions selected by a teacher heuristic, and involves tunable parameters that control the amount of deviation that is allowed. Our approach uses the heuristic as a starting point for further improvements, without the need for any parameters.
\section{The dynamic traveling multi-maintainer problem with alerts}\label{sec:model}
The $K$-DTMPA is a discrete-time model in which a central decision-maker is responsible for selecting the actions of the $K$ maintenance engineers, denoted by $\mathcal{K} = \{1,\ldots, K\}$. To prevent and resolve failures, the assets require maintenance regularly. The engineers maintain a set of assets (machines), denoted by $\mathcal{M} = \{1,\ldots, M\}$, each positioned at a unique location in the network. Each asset $m \in \mathcal{M}$ is subject to degradation which occurs randomly over time. The degradation state is collected in real-time via sensors. After an increase in degradation, an \emph{alert} is issued that informs the decision-maker about the state of the assets. In each time period, the decision-maker selects an action $u_k$ for the $k$-th maintenance engineer, for each $k \in \mathcal{K}$. The action space per engineer consists of actions to travel to another location, idling/continuing or to start maintenance at their location. 

The objective is to minimize the total expected discounted cost over an infinite horizon. Maintenance is referred to as preventive maintenance (PM) when carried out before a failure occurs and it restores the degradation state of the asset to as-good-as-new, as opposed to corrective maintenance (CM) which can be carried out after a failure has occurred and it also restores the state to as-good-as-new. Typically, CM comes at a higher cost compared to PM. Until maintenance is completed, the machine is down during which the decision-maker incurs downtime costs. 
\begin{wrapfigure}{r}{0.5\textwidth}
\centering 
\includegraphics[width=0.85\linewidth]{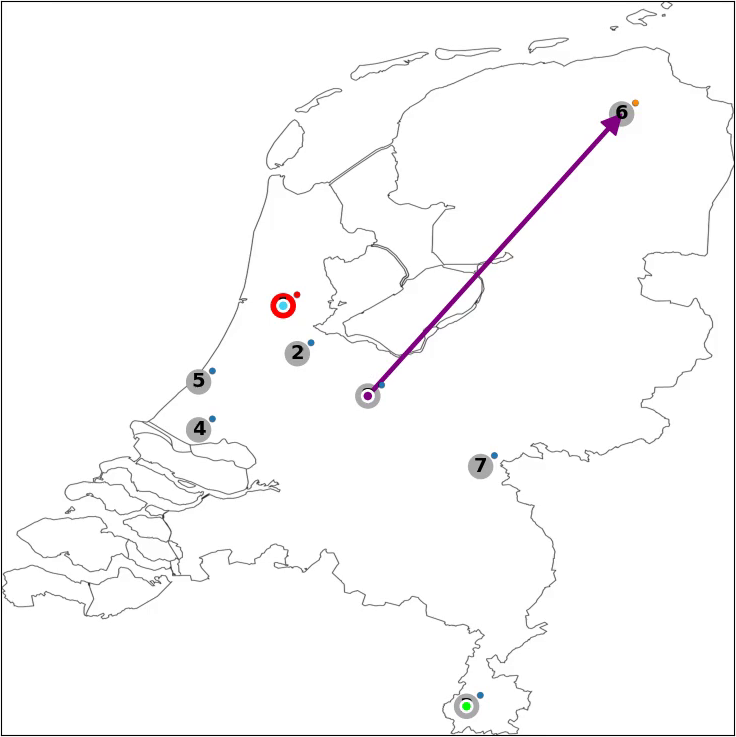}
 \caption{(Figure best viewed in color.) Visualization of the $K$-DTMPA model for an asset network of $M=8$ machines serviced by $K=3$ maintenance engineer. \textit{Blue} dots on top of machine nodes indicate that the machine is healthy, \textit{orange} when alerted or \textit{red} when the machine is down. The engineers are colored \textit{cyan}, \textit{green} and \textit{purple} and are located at Amsterdam, Maastricht and Utrecht, respectively. At discrete decision epochs, engineers can either: (i) idle/continue, (ii) travel to another location or (iii) start maintenance at the current location.}
 \label{fig:k-dtmpa}
\end{wrapfigure}
The connectivity of the network is captured by a travel time matrix $\Theta$ with elements $\theta_{ij}\in\mathbb{N}$ where $i,j \in \mathcal{M}$. Traveling between assets $i$ and $j$ takes $\theta_{ij}$ time units, and requires proportional costs. Each asset $m \in \mathcal{M}$ is assumed to degrade independently which is justified since assets are positioned at different locations. At each time $t\in \mathbb{N}_0$, the asset's degradation state is denoted by $x_m(t)\in \mathcal{N}_m = \{1, \ldots, |\mathcal{N}_m|\}$. Degradation state $x^\textbf{\textrm{h}}_m=1$ represents an as-good-as-new asset and degradation state $x^\textbf{\textrm{f}}_m=|\mathcal{N}_m|$ represents an unavailable asset. When an asset $m$ transitions to the degradation state $x^\textbf{\textrm{f}}_m$, it remains in this degradation state until CM is completed by one of the engineers. In essence, the $K$-DTMPA arises as the multi-maintainer variant of the DTMPA under full (state) information, i.e., the decision-makers are continuously aware of the degradation state of every asset (cf. \citet[Section 3.2, under information level $\textbf{L}_3$]{dtmpa}). See \mbox{Figure \ref{fig:k-dtmpa}} for a visualization of a $K$-DTMPA instance.

The work covered in this paper extends to any degradation model, however, for numerical tractability we have opted to implement the framework propose by \citet{derman1963optimal}: In the numerical section, we assume that, without interference by an engineer, the next state is a ``worse'' state, viz. $x_m(t+1) \geq x_m(t)$ and that the transition time $T^{x_m}_m$ between degradation states $x_m\in \mathcal{N}_m	\setminus\{x^\textbf{\textrm{f}}_m\}$ and $x_m+1$ is random but positive and integer.I.e., the random variable $T^{x_m}_m$ follows a Geometric distribution with success parameter $p_m^{x_m} \in (0,1]$. 

The remainder of this section formalizes the sequential decision-making $K$-DTMPA model framework. We detail the states, actions, transitions, costs and the optimization objective hereinafter. 

\subsection{States, actions and transitions}\label{subsec:states-actions-transitions}
The \textit{state} $h\in \mathcal{H}$ contains the degradation state of all machines in the network and $K$ blocks capturing the status of the engineers. Each block consists of three elements that describe the engineer's location, current activity and availability. A state can thus be represented as a vector $h = (x_1, \ldots, x_M, \ell_1, \iota_1, \delta_1, \ldots, \ell_K, \iota_K, \delta_K) \in \mathcal{H}$, with a minor abuse of notation. Here, $x_m\in\mathcal{N}_m$ represents the degradation state of asset $m\in\mathcal{M}$; $\ell_k\in\mathcal{M}$, $k\in\mathcal{K}$, denotes the location of the $k$-th engineer; $\iota_k \in I = \{0,1\}$ indicates whether this engineer is currently carrying out maintenance (otherwise, the engineer is either traveling or idling); $\delta_k \in\Delta \subset \mathbb{N}_0$ counts the remaining time units that the engineer is occupied (e.g., $\delta_k = 0$ specifies that the engineer is available). Under the assumption of Geometric transition times between subsequent degradation states, the elapsed times in the state space description $(h)$ can be excluded.

At every time instance, for each maintenance engineer, given $h \in \mathcal{H}$, the decision-maker must choose one of the following options: (i) start traveling to location $m\in\mathcal{M}\setminus\{\ell_k\}$, say $u_m$; (ii) start maintenance at the present location, say $v$; or (iii) continue the ongoing activity or remain idle, say $u_{\ell_k}$, with $\ell_{k}$ denoting the location of the engineer for which an action is selected, say the $k$-th engineer, $k \in \mathcal{K}$. In particular, when $\delta_k > 0$, the $k$-th engineer is unavailable and therefore action $u_{\ell_k}$ must be chosen, while if $\delta_k=0$, then the action can be chosen from the set $\{u_m\}_{m\in \mathcal{M}} \cup \{v\}$, where maintenance action $v$ is only available if the machine at its location is not already being maintained by another engineer. Thus, the \textit{state-dependent action set for the $k$-th maintenance engineer} becomes
$$\mathcal{U}_k(h) = \begin{cases}
\{u_m\}_{m\in \mathcal{M}} \cup \{v\} & \textrm{if } \delta_k = 0 \land \sum_{k'\in\mathcal{K}\backslash\{k\}} \mathds{1}_{\{\ell_{k'} = \ell_k, \iota_{k'} = 1 \}} = 0, \,\\
\{u_m\}_{m\in \mathcal{M}} & \textrm{if } \delta_k = 0 \land \sum_{k'\in\mathcal{K}\backslash\{k\}} \mathds{1}_{\{\ell_{k'} = \ell_k, \iota_{k'} = 1 \}} > 0, \,\\
\{u_{\ell_k}\} & \textrm{if } \delta_k > 0.
\end{cases}$$

\emph{The state-dependent action set} $\mathcal{U}(h)$ is the Cartesian product of the $K$ individual state-dependent action sets excluding those actions that result in two or more engineers simultaneously maintaining a machine, i.e.,

\begin{equation*}
 \mathcal{U}(h) = \Big( \mathcal{U}_1(h) \times \ldots \times \mathcal{U}_K(h) \Big) \setminus \{a \in \mathcal{U}_1(h) \times \ldots \times \mathcal{U}_K(h) \mid \sum_{ k \neq k' } \mathds{1}_{ \{a_k = v, a_{k'} = v, \ell_k = \ell_{k'} \}} > 0 \}.
\end{equation*}
The state transition $h \rightarrow h'$ is decomposed in two stages: The first stage is determined by the deterministic consequences of the chosen actions, say $a = (a_1, \ldots, a_K) \in \mathcal{U}(h)$, and the second stage is determined by the random evolution of the degradation processes. More specifically, $h \rightarrow h'$ is decomposed into $h \overset{a_1}{\rightarrow} h^{a_1} \overset{a_2}{\rightarrow} \ldots \overset{a_K}{\rightarrow} h^{a_K} =: h^a$ and to $h^a \overset{t \rightarrow t+1}\longrightarrow h'$. The order of handling the actions $a_1, \ldots, a_K$ is irrelevant, therefore, we only detail how the action of the $k$-th maintenance engineer is processed.

In the case that $(a_k = u_{\ell_k})$, the engineer remains at the current location $\ell_k$ and continues the ongoing action or idles. The remaining unavailability $\delta_k$ is increased by one if the engineer is idle. When $(a_k = u_m \land m \neq \ell_k)$, the engineer starts to travel to location $m$. The remaining unavailability $\delta_k$ is increased with the travel time $\theta_{\ell_k m}$ and the location $\ell_k$ is updated accordingly. The third and fourth possibility, when $(a_k = v \land x_{\ell_k} \neq x^\textbf{\textrm{f}}_{\ell_k})$ or $(a_k = v \land x_{\ell_k} = x^\textbf{\textrm{f}}_{\ell_k})$ respectively, represent initiating PM or CM at the current location, depending on the status $x_{\ell_k}$ of the machine at the location of the engineer. The remaining unavailability $\delta_k$ is increased with the duration of maintenance, which is either $t^{\textrm{PM}}_{\ell_k} \in \mathbb{N}^+$ or $t^{\textrm{CM}}_{\ell_k} \in \mathbb{N}^+$. Following the modeling assumptions, the degradation state $x_m$ of machine $m$ is set to $x^\textbf{\textrm{f}}_m$ to indicate that the machine is unavailable during maintenance.

For $h^a = (x^a_1, \ldots, x^a_M, \ell^a_1, \iota^a_1, \delta^a_1, \ldots, \ell^a_K, \iota^a_K, \delta^a_K)
\overset{t \rightarrow t+1}\longrightarrow h' = (x'_1, \ldots, x'_M, \ell'_1, \iota'_1, \delta'_1, \ldots, \ell'_K, \iota'_K, \delta'_K )$, determining $h'$ is separated in two steps: First, we update the state of every machine according to the random evolution of the degradation process or completion of maintenance. Subsequently, we update the remaining state variables. In more detail, the state is modified as follows:
\begin{enumerate}
\item One of the following triptych of cases determines the evolution of machine $m\in\mathcal{M}$:
If $(\ell^a_k, \iota^a_k, \delta^a_k) = (m, 1, 1)$ for some $k\in\mathcal{K}$, then $x_m'=x^\textbf{\textrm{h}}_m$. This first case represents that after completion of PM or CM, the state of the machine is updated to as-good-as-new. Else, with probability $\mathbb{P}(T_m^{x^a_m} = 1)$, $x_m'= x^a_m+1 $, the machine transitions to the subsequent degradation state. This is the second case. Otherwise, in the third case, which occurs with probability $1-\mathbb{P}(T_m^{x^a_m} = 1)$, $x_m' = x^a_m $, the machine degradation state remains the same.
\item The evolution of each engineer $k\in\mathcal{K}$ is deterministic: $(\ell^a_k, \iota^a_k\mathds{1}_{\{ \delta^a_k > 1\}}, \delta^a_k - 1) \overset{t \rightarrow t+1}\longrightarrow (\ell_k', \iota_k', \delta_k')$; i.e., the remaining unavailability $\delta_k'$ of the $k$-th engineer is decreased by one when continuing an ongoing activity. The indicator $\iota_k'$ resets when the $k$-th engineer completes an activity.
\end{enumerate}

\subsection{Cost structure and objective}\label{subsec:costsobjective}
The cost structure of the $K$-DTMPA includes costs for travel, maintenance and asset unavailability. A small cost $c^{\textrm{T}} \in\mathbb{R}^+$ is paid for each unit of time that an engineer travels, independent of the maintenance engineer, the origin and the destination. Initiating PM or CM on machine $m\in\mathcal{M}$ costs $c^{\textrm{PM}}_{m} \in\mathbb{R}^+$ or $c^{\textrm{CM}}_{m} \in\mathbb{R}^+$, respectively. CM is assumed to be more costly, viz. $c^{\textrm{CM}}_{m} \geq c^{\textrm{PM}}_{m}$. The failed state $x^\textbf{\textrm{f}}_m$ models the unavailability of asset $m$, i.e., when the asset has failed or is under repair. The cost of downtime is $c^{\textrm{DT}}_{m}\in\mathbb{R}^+$ per time unit, regardless of the source of disruption. Thus, when taking action $a \in \mathcal{U}(h)$ in state $h$, the incurred costs are: 
\begin{align*}C(h,a) &= \sum_{k\in\mathcal{K}} \left(c_{\ell_k}^\textrm{CM}\mathds{1}_{\{a_k = v, x_{\ell_k} = x^\textbf{f}_{\ell_k}\}} + c_{\ell_k}^\textrm{PM}\mathds{1}_{\{a_k = v, x_{\ell_k} < x^\textbf{f}_{\ell_k}\}} + c^\textrm{T} \left(\mathds{1}_{\{a_k \neq \ell_k, a_k \neq v\}} + \mathds{1}_{\{\delta_k > 0, \iota_k = 0\}} \right)\right) \\ &+ \sum_{m\in\mathcal{M}} c^\textrm{DT}_m \mathds{1}_{\{x_m = x^\textbf{f}_m \}} .\end{align*}
The objective is to devise a policy $\pi$ that minimizes the total expected discounted cost. A policy $\pi = (\pi_1, \pi_2, \ldots,\pi_t, \ldots )$ is defined as a sequence of decision rules, where the decision rule $\pi_t$ is a probability distribution over the action space $\mathcal{U}(h)$ at time $t$, given the state $h \in \mathcal{H}$. We denote with $\pi^k_t$ the induced probability distribution over the action set for the $k$-th engineer $\mathcal{U}_k(h)$. Let $\gamma \in [0, 1)$ be the discount factor and $J(\pi)$ be the total expected discounted cost. Thus, the objective is to determine the optimal policy $\pi_*$ satisfying \begin{equation*}
\label{eq: objective}
 \pi_* = \textrm{arg }\underset{\pi}{\textrm{min }} J(\pi) = \textrm{arg }\underset{\pi}{\textrm{min }} \lim_{T\to\infty} \mathbb{E}_{\pi} \Bigg[ \sum_{t=0}^T \gamma^t C\left( h(t), a(t) \right) \Bigg| h(0)=h\Bigg],
\end{equation*}
where $(h(t),a(t))$ denotes the tuple of the state and the
respective action given the policy $\pi_t$ at time $t$, $t \geq 0$, and $C(\cdot)$ denotes the associated cost (maintenance, travel and downtime).

\section{Heuristic policies and benchmarks}\label{sec:heuristics}
In this section, we discuss the aspects that characterize a good policy for a $K$-DTMPA instance. Subsequently, we detail the dispatching heuristic with which we equip existing ranking heuristics. We argue that the resulting class of heuristics contains compelling benchmarks for suitable subclasses of $K$-DTMPA instances. The $K$-DTMPA framework encompasses a wide range of dynamic traveling maintainer problem instances, however, suitable parameter choices will reduce the resulting problem to well-studied problems. For instance, by considering a single maintainer, the problem reduces to the DTMPA under full state information studied by \citet{dtmpa}. Setting the number of states $N_m \equiv 2$ yields a problem similar to a dynamic vehicle routing problem. Thus, we can compare our DRL algorithm to heuristics designed for such special cases. 

\subsection{Aspects of a good policy}\label{subsec:policy_aspects}
What constitutes a good policy is the ability to jointly consider: the network layout, the spatial and temporal information of all engineers, the uncertainty in the evolution of each machine's condition, and the cost structures. Such a policy must account for the fact that faster degrading machines typically require maintenance more regularly and expensive CM actions (compared to PM) must be avoided. In addition, the downtime costs require engineers to move proactively to anticipate future events in the network.

\begin{figure}[!ht]
\centering
 \begin{subfigure}[b]{0.3\columnwidth}
 \includegraphics[width=\linewidth]{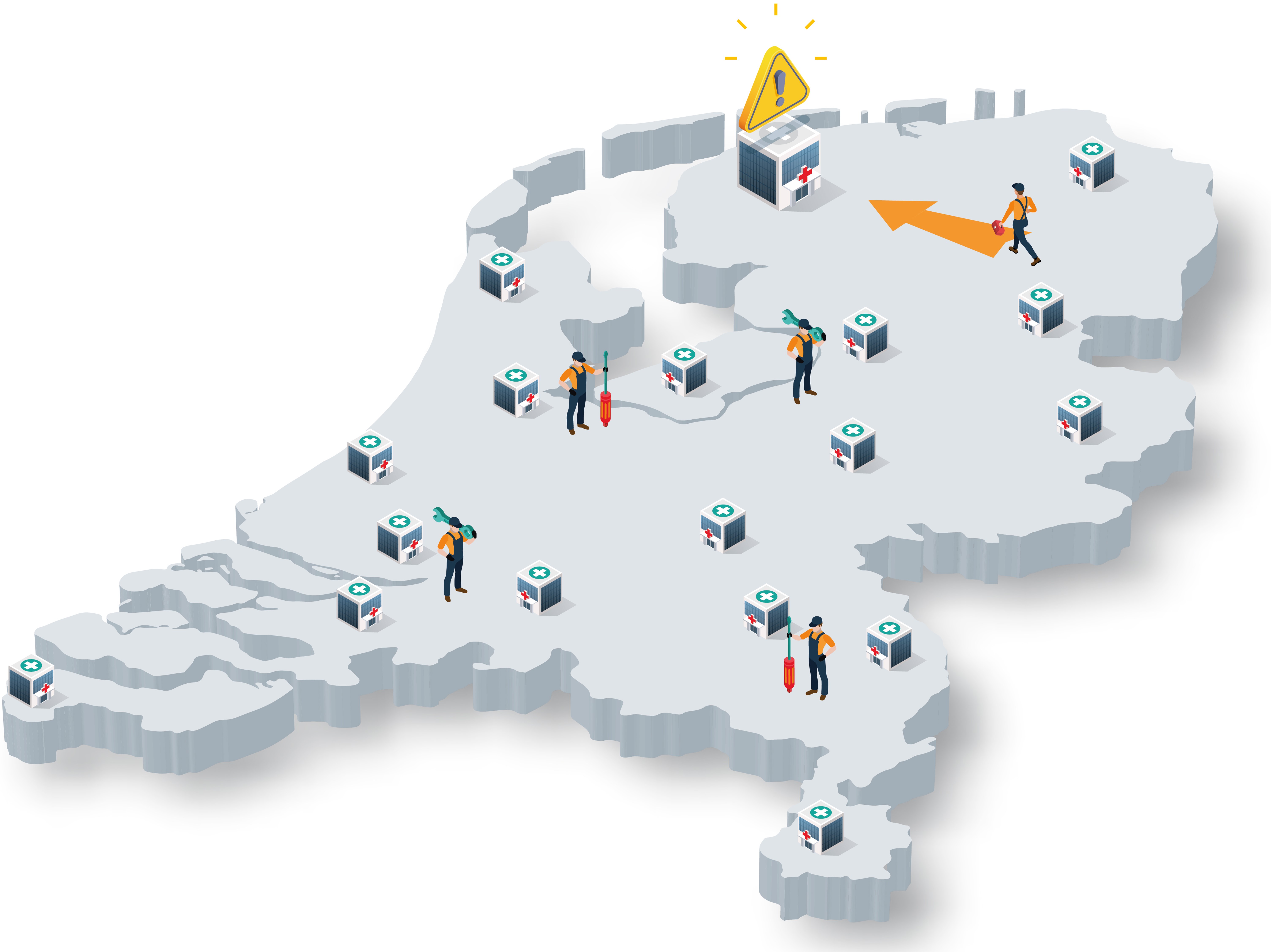}
 \caption{Failure: When assets fail, \\available engineers must be dispatched efficiently. }
 \label{fig:2}
 \end{subfigure}
 \hfill %%
 \begin{subfigure}[b]{0.3\columnwidth}
 \includegraphics[width=\linewidth]{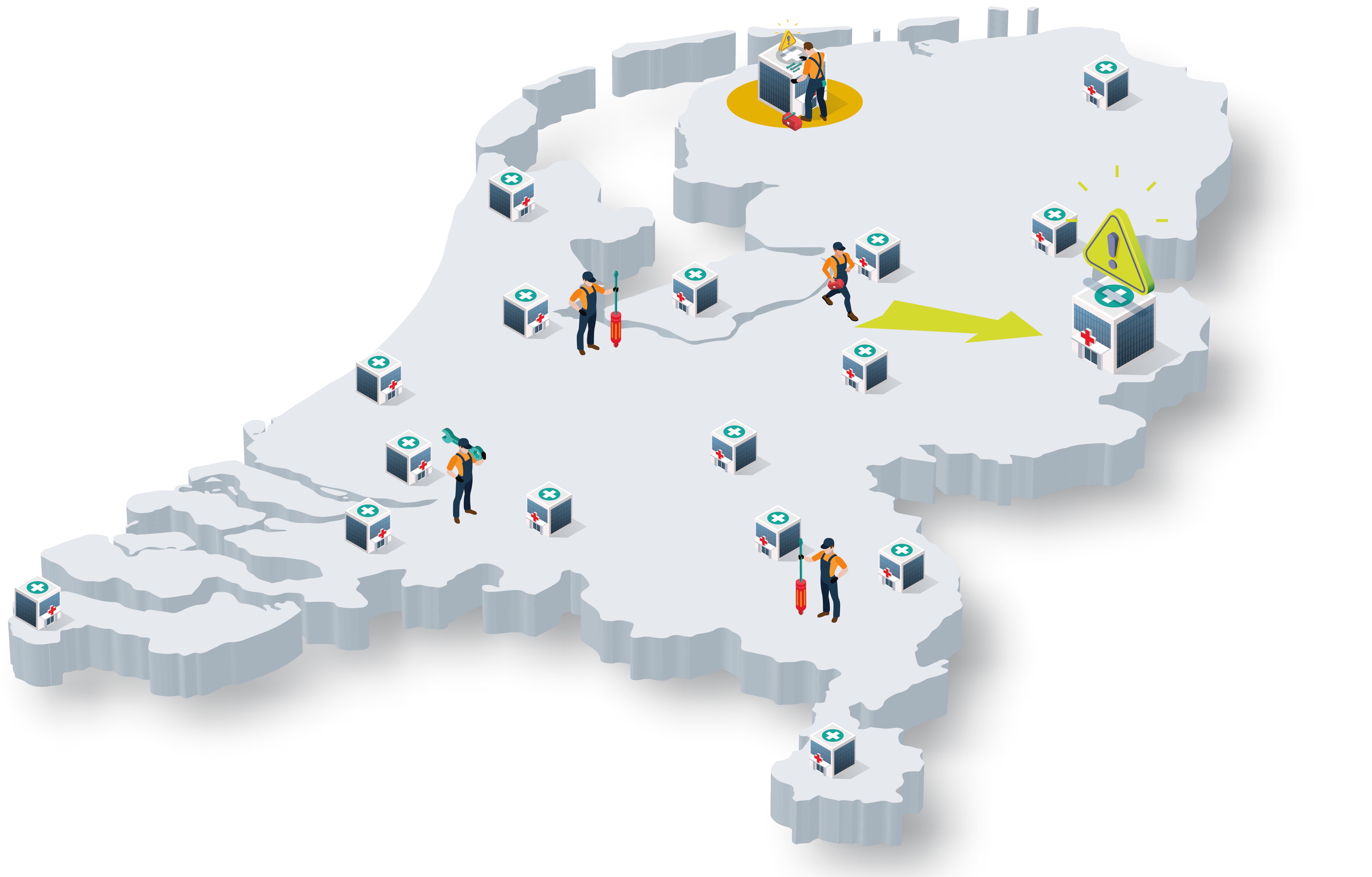}
 \caption{Alert: When an alert is issued, the decision-maker must conduct a risk \\urgency assessment to decide whether to dispatch an engineer.}
 \label{fig:3}
 \end{subfigure}
 \hfill
 \begin{subfigure}[b]{0.3\columnwidth}
 \includegraphics[width=\linewidth]{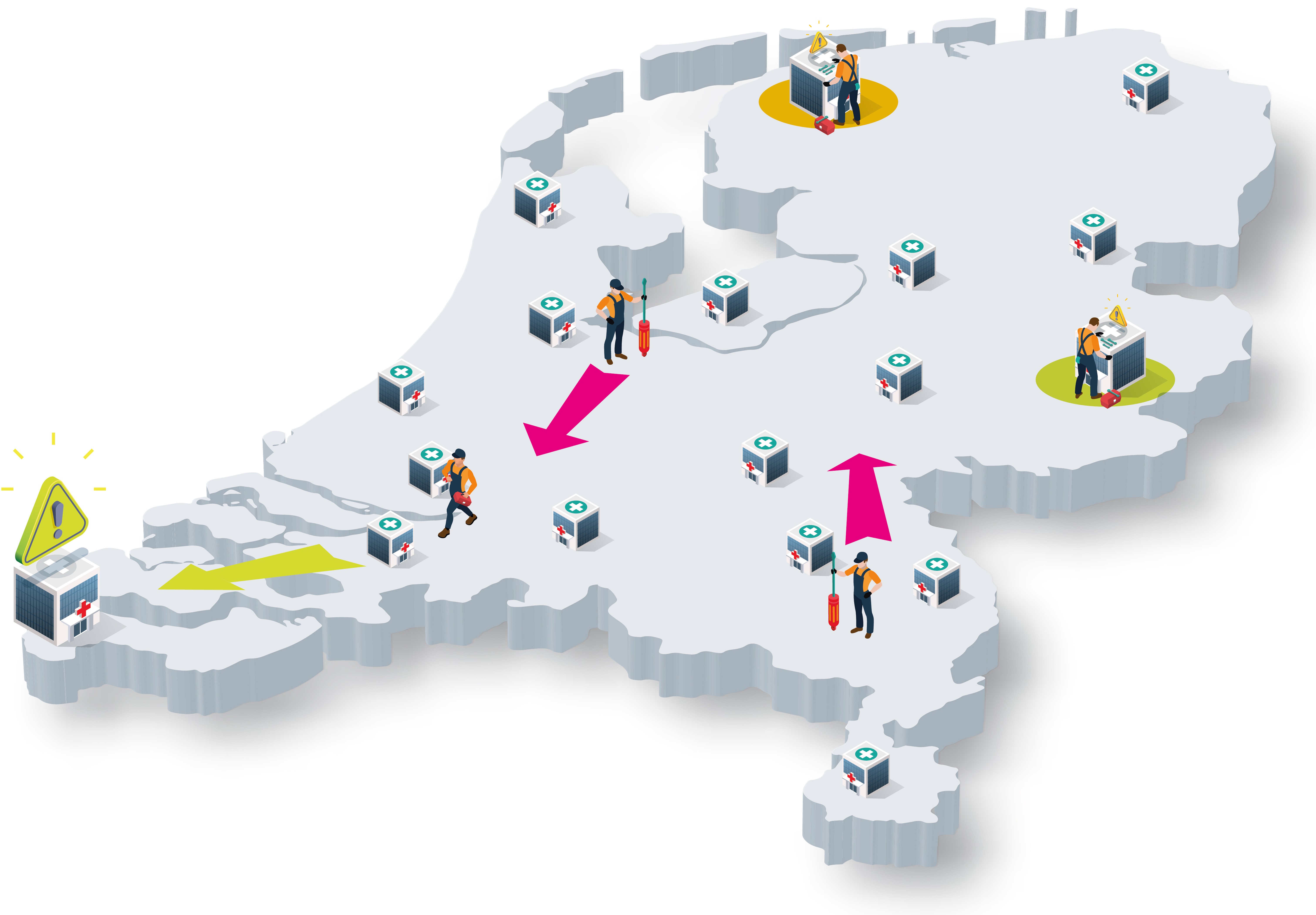}
 \caption{Repositioning: Idle engineers are proactively repositioned to be closer to future alerts and failures.}
 \label{fig:4}
 \end{subfigure}
\caption{Visualization and description of the three envisioned policy aspects.}
\label{fig:policy_aspects}
\end{figure}

On that account, the aspects that construct a good policy are envisioned to be the following: (i) Efficient dispatching of the engineers to the various locations; (ii) Assessing the risk of delaying preventive maintenance; and (iii) Tactical repositioning of any remaining available engineers. See Figure \ref{fig:policy_aspects} for the visualization and further elaboration on these three policy aspects.

\subsection{Dispatching heuristic policy}\label{subsec:greedy_heur}
The existing greedy and reactive ranking heuristic solutions introduced by \citet[Section 5.1]{dtmpa} are not immediately applicable to the $K$-DTMPA setting since they do not include a dispatching mechanism and since they operate at an information disadvantage. In the setting of \citet{dtmpa}, these heuristics only observe transitions to the first degradation state $x_m^\textbf{h}+1$ and the failed state $x_m^\textbf{f}$, while in the $K$-DTMPA setting all degradation state transitions are observed. Therefore, it is paramount to extend the greedy and reactive ranking heuristics to the multi-maintainer setting cf. the policy aspects discussed in \mbox{Section \ref{subsec:policy_aspects}}.
The \emph{reactive} heuristic maintains only failed machines whereas the \emph{greedy} heuristic maintains both alerted and failed machines. To determine which machines to maintain, we require a \emph{ranking} between assets. In this work, we consider a state-dependent threshold ranking leveraging the available degradation information. A second optimization step accounts for costs and travel times when dispatching the available maintenance engineers to the ranked assets. The remaining engineers are assigned the action to remain idle or continue with the current activity, i.e., no repositioning step is performed. We denote the resulting dispatching heuristic policy by $\pi_\textrm{D}$.

{\em State-dependent threshold ranking.} Given the alert information, the machines are ranked on their observed degradation level. Let $t\in\mathbb{N}_0$ be the current time and $h \in \mathcal{H}$ be the state at time $t$. Machine $m$ is added to the ranking when $x_m(t) \geq s_m(h)$, where $s_m(h) \in \mathcal{N}_m$ is the state-dependent degradation threshold corresponding to machine $m$. Note that setting $s_m(h) \equiv |\mathcal{N}_m|$ yields a reactive policy. If there are more ranked assets than $K'$ available engineers, we iteratively reduce the ranking as follows: An asset that is farthest from one of the available engineers is chosen randomly and it is removed from the ranking.

{\em Dispatching step.} The ranked assets must now be assigned to the available engineers. The formulation of the assignment problem in its general form is as follows: The problem instance has a number of engineers and a number of maintenance jobs. Any engineer can be assigned to perform any job, incurring some cost that may vary depending on the engineer-job assignment. It is required to perform as many jobs as possible by assigning at most one engineer to each job and at most one job to each engineer in such a way that the total cost of the assignment is minimized. By construction, the assignment problem contains at most $K'$ jobs. When there are fewer jobs than available engineers, the so-called unbalanced assignment problem can be reformulated as a balanced assignment problem by adding dummy jobs (i.e., jobs with a cost of 0 for each available engineer). The pairwise travel time between ranked assets and available engineers seems to be the rational choice to construct the cost matrix in the case of identical assets (in terms of cost parameters and distributional degradation characteristics). The Hungarian method solves the constructed assignment problem in $\mathcal{O}(M^3)$ polynomial time complexity \citep[Chapter 17.2]{schrijver2003combinatorial}. The engineers are dispatched according to the solution to the constructed assignment problem.

\section{Approximate policy iteration for $K$-DTMPA}\label{sec:api}
\Citet{dtmpa} have shown that DRL, specifically, $n$-step quantile regression double Q-Learning ($n$QR-DDQN), produces near-optimal policies for DTMPA instances. However, training times exceed 12 hours even for networks containing up to 6 assets. To learn policies for $K$-DTMPA instances, we adopt a form of approximate policy iteration (API), which enables us to distribute the sample collection over multiple compute nodes (cf. \cite{chess}), contributing substantially to scalability. In particular, we adopt deep controlled learning (DCL) \citep{temizöz2023deep}, which combines variance reduction and optimized allocation of roll-outs for greater efficiency and trains neural networks using cross-entropy (instead of Euclidean) loss functions. DCL has been shown to outperform other DRL algorithms such as \emph{proximal policy optimization} or \emph{asynchronous advantage actor-critic} on inventory problems \citep{temizöz2023deep}.

To successfully apply API/DCL to $K$-DTMPA instances, we combine several novel ideas that shall be discussed in depth in this section. First, we provide an overview of the application of API to $K$-DTMPA instances (for details, we refer to \cite{temizöz2023deep}) and then discuss the feature representation of state information, training the neural network classifier and suitable initial solutions in the forthcoming sections.

\begin{figure}[ht]
\centering \includegraphics[width=0.7\linewidth]{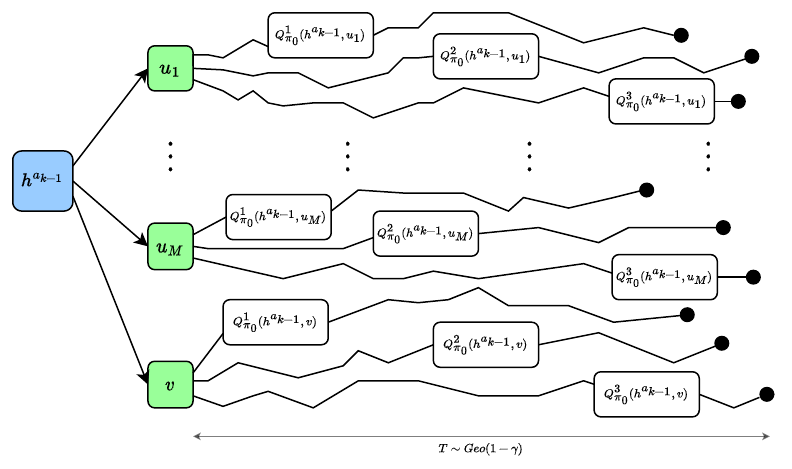}
 \caption{Visualization of the simulation-based policy $\hat{\pi}^+$. In state $h^{a_{k-1}}$, the policy prescribes to follow the action $\tilde{a} \in \mathcal{U}_k(h^{a_{k-1}})$ (recall that $\mathcal{U}_k(h^{a_{k-1}}) \subseteq \{u_{m}\}_{m=1}^M \cup \{v\})$ that minimizes the average undiscounted trajectory cost $\hat{q}_{\pi_0}(h^{ a_{k-1}}, \tilde{a})$. Unbiased estimates of the action-value function are computed from $r$ independent roll-out simulations whose length follows a Geometric distribution with parameter $1-\gamma$.}
 \label{fig:api4kdtmpa}
\end{figure}

Recall from \mbox{Section \ref{subsec:states-actions-transitions}} that the action space $\mathcal{U}(h)$ grows exponentially in the number of engineers $K$. To vastly reduce the action space complexity, we train a neural network to select the actions for the engineers sequentially in a fixed order. Due to symmetry, any ordering of the engineers can be adopted. To enable cooperation, after each action selection, the input is updated with the consequences of the action (see \mbox{Section \ref{subsec:states-actions-transitions}}) before selecting an action for the next engineer. Let $\mathcal{I}= \{1,\ldots,|\mathcal{I}|\}$ be some index set. Starting from some initial policy $\pi_0$, we interact with the environment to collect a data set $\mathcal{D} = \{ (h^{ a_i^{k-1}}_i, a^k_i) \mid i\in \mathcal{I}\}$ consisting of state-action tuples $(h^{ a_i^{k-1}}_i, a^k_i)$ for which the action has been obtained using a form of \emph{simulation optimization}. More specifically, starting from state $h_i$, given that the first $k-1$ maintenance engineers select actions $a^1_i, \ldots, a^{k-1}_i$, we select the action $a^k_i$ for the $k$-th engineer in state $h_i^{ a^{k-1}_i}$ that minimizes the action-value function $q_{\pi_0}(h_i^{ a^{k-1}_i},\cdot)$ associated with policy $\pi_0$. For $\tilde{a} \in \mathcal{U}_k(h_i^{a_i^{k-1}})$, the action-value function is defined as follows:
\begin{equation*}
 q_{\pi_0}(h_i^{ a^{k-1}_i},\tilde{a}) = \lim_{T\to\infty} \mathbb{E}_{\pi_0} \Bigg[ \sum_{t=0}^T \gamma^t C\left( h(t), a(t) \right) \Bigg| h(0)=h_i, a_1(0) = a^1_i, \ldots, a_{k-1}(0) = a^{k-1}_i, a_k(0) = \tilde{a} \Bigg].
\end{equation*}
In other words, $q_{\pi_0}(h_i^{ a^{k-1}_i},\tilde{a})$ is the total expected discounted cost when selecting action $\tilde{a}$ for state $h_i^{ a^{k-1}_i}$ and following the policy $\pi_0$ in the remainder of the roll-out, i.e., $a_{j}(0) \sim \pi^j_0(h_i^{a_{j-1}(0)})$ for $j=k+1,\ldots,K$ and $a(t) \sim \pi_0(h(t))$ for all $t\geq1$. To estimate $q_{\pi_0}(h_i^{ a^{k-1}_i},\tilde{a})$, we generate $r \in \{r_\textrm{min}, \ldots, r_\textrm{max}\}$ independent roll-out simulations of length $T\sim \textrm{Geo}(1-\gamma)$ and compute the undiscounted trajectory costs $Q^j_{\pi_0}(h_i^{ a^{k-1}_i}, \tilde{a})$ for $j=1,\ldots,r$. The resulting unbiased estimator $\hat{q}_{\pi_0}(h_i^{ a^{k-1}_i}, \tilde{a})$ \citep{HAVIV1992267} satisfies
\begin{equation*}
\hat{q}_{\pi_0}(h_i^{ a^{k-1}_i}, \tilde{a}) = \frac{1}{r} \sum_{j = 1}^{r} Q^j_{\pi_0}(h_i^{ a^{k-1}_i}, \tilde{a}).
\end{equation*}
Subsequently, the improved action $a^k_i$ for state $h_i^{ a^{k-1}_i}$ is determined as $$a^k_i = \underset{\tilde{a} \in \mathcal{U}_k(h_i^{ a^{k-1}_i})}{\textrm{arg min}} \hat{q}_{\pi_0}(h_i^{ a^{k-1}_i}, \tilde{a}) =: \hat{\pi}^+(h_i^{ a^{k-1}_i}).$$
We refer to $\hat{\pi}^+(h_i^{ a^{k-1}_i})$ as the simulation-based policy for state $h_i^{ a^{k-1}_i}$, see \mbox{Figure \ref{fig:api4kdtmpa}} for a visualization.
Ideally, we include those states that are visited most frequently under this simulation-based policy. Thus, starting from some initial state $h_0 \in \mathcal{H}$, states for $\mathcal{D}$ are those that are encountered when selecting in each decision epoch, for each engineer, the randomized simulation-based policy: The policy that chooses a random action with probability $\epsilon \in [0,1]$ and with probability $1 - \epsilon$, follows the simulation-based policy $\hat{\pi}^+$.

Subsequently, we train a neural network classifier on $\mathcal{D}$ (see \mbox{Section \ref{subsec:training_NN}} for more details) which induces a hopefully improved policy. API can be transformed into an iterative scheme by collecting new data using the improved neural network policy. Like exact policy iteration, API can potentially improve any heuristic policy and find good solutions in a handful of iterations. 

Summarized, the API algorithm consists of the following three steps:
\begin{enumerate}
\item Choose a suitable initial solution $\pi_0$.
\item Construct the data set $\mathcal{D}$ using $\pi_0$.
\item Train a neural network classifier on the constructed data set $\mathcal{D}$.
\end{enumerate}
For step three above, the neural network can be interpreted as a parameterized function from $\mathbb{R}^m$ to $\mathbb{R}^n$ for some $m,n\in\mathbb{N}$. Let such a (generic) function be denoted by $N_\theta(\cdot)$, where $\theta$ denotes the function parameters. The input to the neural network is the feature representation $f(h)\in \mathbb{R}^m$ of a state $h\in\mathcal{H}$. The output of the neural network $N_\theta(\cdot) \in \mathbb{R}^{n}$, where $n = M+1$, is transformed into a probability distribution over the action space and the action $\tilde{a}$ which is assigned the highest probability $N_\theta( \cdot )_{\tilde{a}}$ is chosen, i.e., the neural network induces a policy. Given the actions $a_1,\ldots,a_{k-1}$ for the first $k-1$ engineers, the action $a_k$ for the $k$-th engineer in state $h$ is thus determined from the following decision rule: 
\begin{equation*}
 \pi^k_{\theta}( f(h^{ a_{k-1}}) ) = \underset{\tilde{a}\in \mathcal{U}_k(h^{ a_{k-1} })}{\textrm{arg max}} [(N_\theta( f(h^{ a_{k-1} }) )_{\tilde{a}}],
\end{equation*}
where by convention $h^{ a_{0}} = h$. We denote by $\pi_{\theta}$ the neural network policy that selects in every decision epoch, for each maintenance engineer $k$, the action $\pi^k_{\theta} (f(h^{ \pi^{k-1}_{\theta}}))$. 

\subsection{Feature representation}\label{subsec:feature_representation}
We propose a handcrafted feature representation to make the state information suitable for input into the neural network. Although a much more compact representation of the state $h$ is possible, we propose a feature design based on conveying the state information per location in a form that is \emph{tailored to the engineer for which we are currently selecting an action}. We have found that this engineer-centric feature design is crucial to efficiently learn cooperative dispatching mechanisms. Specifically, the state is transformed as follows:
\begin{align*}
f_1(h) = \left( x_1, n^\textrm{av}_1, n^\textrm{ua}_1, t^\nu_1, t^{\Theta_1}_1, t^{\Theta_2}_1, \xi_1, \ldots, x_M, n^\textrm{av}_M, n^\textrm{ua}_M, t^\nu_M, t^{\Theta_1}_M, t^{\Theta_2}_M, \xi_M, \sum_{m\in\mathcal{M}} n^\textrm{av}_m \right),
\end{align*}
i.e., the feature vector $f_1(h)$ contains an information block $(x_m, n^\textrm{av}_m, n^\textrm{ua}_m, t^\nu_m, t^{\Theta_1}_m, t^{\Theta_2}_m)$ that can be computed from $h$ for each $m\in\mathcal{M}$, and one additional feature. Here, $x_m$ is the observed degradation level of machine $m$, $n^\textrm{av}_m$ and $n^\textrm{ua}_m$ denote the number of available and unavailable maintenance engineers at location $m$, respectively. The entry $t^\nu_m$ captures the remaining time to the completion of a maintenance job at location $m$, whilst $t^{\Theta_1}_m$ and $t^{\Theta_2}_m$ indicate the remaining travel time until the first and second arrival of an engineer at location $m$, respectively. By convention, the default value of $t^\nu_m, t^{\Theta_1}_m$ and $t^{\Theta_2}_m$ is $0$. The last block entry $\xi_m$ indicates whether the maintenance engineer for which we are currently selecting an action is present at location $m$. Lastly, the total number of available engineers is added as an additional feature. All in all, the dimension of the feature vector equals $n = 7M+1$ and is independent of the number of maintenance engineers $K$. In \mbox{Section \ref{subsec:feature_design}}, we compare the quality of the trained neural networks using the proposed feature representation against (i) a similar feature representation albeit without the last feature (say $f_2(h)$), i.e., $f_1(h) = \left( f_2(h), \sum_{m\in\mathcal{M}} n^\textrm{av}_m \right)$, and (ii) against the most compact state representation (say $f_3(h)$).

\subsection{Training the neural network classifier}\label{subsec:training_NN}
API relies on supervised learning to train neural networks. In our context, supervised learning finds a relation between the actions $\tilde{a} \in \mathcal{U}_k( h^{a_{k-1}} )$ as taken by the simulation-based policy $\hat{\pi}^+$ and the feature representation $f(h^{a_{k-1}})$ of the state $h^{a_{k-1}} \in \mathcal{H}$. Specifically, we employ a multilayer perceptron consisting of $L\in\mathbb{N}$ layers. In each layer $l \in \{1, \ldots, L\}$, an affine transformation of the input is combined with a nonlinear activation function. For our experiments, we adopt a standard supervised learning algorithm cf. \citet[Appendix A]{temizöz2023deep}. For learning the parameters $\theta$, we split the data set $\mathcal{D}$ in a training set and a test set. The training loss $L(\theta)$ measures the ``distance” between $\hat{\pi}^+$ and the neural network policy $\pi_\theta$ for the states in the training set. Fitting $\theta$ is an iterative, gradient-based process: In each step, the gradient of $L(\theta)$ with respect to $\theta$ is estimated, and subsequently, $\theta$ is updated by taking a step in the opposite direction. We terminate when the loss on the test set, defined analogously to the loss for the training set, no longer decreases.

\subsection{Initial solutions}\label{subsec: initial_sol}
Our experiments reveal that choosing a suitable policy $\pi_0$ to initiate API/DCL may significantly reduce computation times. In particular, we have found that the following four properties play a key role.

{\em Exploration.} Sufficiently many states must be encountered to yield a rich data set $\mathcal{D}$, dismissing for instance the idle policy, viz. the policy that always selects the action to idle for every engineer.

{\em Cooperation.} Desirable cooperative behavior is typically hard to learn or improve, limiting the use of the random policy to single-maintainer instances.

{\em Computational complexity.} The ability to generate large data sets in a reasonable time is of paramount importance for solving large-scale $K$-DTMPA instances.

{\em Non-self-correcting.} Policies that revert deviations are typically not suitable. For example, suppose we adopt a network decomposition approach as $\pi_0$. Such an approach assigns engineers to predetermined clusters of machines (more details are provided in \mbox{Section \ref{sec:numerical_experiments}}). Under such a policy, dispatching an engineer outside its cluster will typically be followed by a correcting action that returns the engineer to its cluster, which inhibits the effective learning of cooperative behavior.

The dispatching heuristic policies developed in \mbox{Section \ref{subsec:greedy_heur}} satisfy all four properties. The dispatching heuristics encounter all the relevant states while maintaining the network. Moreover, the cooperative behavior of the engineers is optimized in a myopic fashion in polynomial time. The heuristics are also non-self-correcting since they implement no special asset-maintainer constraints.
\section{ Numerical experiments }\label{sec:numerical_experiments}

To assess the performance and scalability of the algorithm proposed in Section \ref{sec:api}, we construct several asset networks with the number of machines $M$ ranging from $4$ to $35$. Machine degradation times, i.e., times to go from one degradation level to the next, are geometrically distributed. Under these circumstances, the $K$-DTMPA is a large-scale computationally intractable MDP. We assess the performance of DCL on a selection of $K$-DTMPA instances for which compelling benchmarks are available: single maintainer instances and D\&R instances. Indeed, the heuristic approaches developed in \mbox{Section \ref{subsec:greedy_heur}} are specifically suitable for the latter instances. The preventive maintenance instances serve as more complex cases to learn new valuable insights using DCL.

\paragraph{Single maintainer instances}
In the situation that $K=1$, the $1$-DTMPA reduces to a DTMPA under full state information, cf. \citet[Section 3.2]{dtmpa}. For small instances containing up to $M=4$ machines, the optimal policy is available as a benchmark. 

\paragraph{Dispatching \& repositioning instances}
D\&R instances are $K$-DTMPA instances on networks where all assets are identical, both in terms of cost structure and degradation dynamics. Moreover, all assets are assumed to have only two states: healthy and failed. As such, the element of preventive maintenance is eliminated and the sole objective becomes to minimize the unavailability of the assets. The dispatching heuristic developed in \mbox{Section \ref{subsec:greedy_heur}} has exactly this objective in mind and can be improved by selecting additional repositioning actions. Therefore, the dispatching heuristic will serve both as the benchmark and as the initial policy for DCL. 

\paragraph{Preventive maintenance instances}
We modify the D\&R instances by including an additional state, in total we have three states: healthy, degraded and failed. The machine will transition to the degraded state on average at 75\% of the machine's life expectancy. This complicates the objective as the policy now needs to jointly consider the cost structures and the network layout, including the position and availability of the engineers. For such instances, no strong benchmark exists in prior work. As such, we propose a traditional heuristic approach by means of network decomposition: The $K$-DTMPA instance is decomposed into $K$ disjoint $1$-DTMPA instances, which can each individually be optimized using DCL, which is known to perform well for $1$-DTMPA instances. We shall refer to this policy as $\pi^\textrm{DEC}_\textrm{D}$. To compose the $i$-th generation policy $\pi^\textrm{DEC}_{\theta_i}$, per cluster, we select the best-found neural network policy so far. An example of a decomposition of a $K$-DTMPA instance is given in \mbox{Section \ref{subsec: networks}}.

A detailed setup of the experiments follows in the remainder of this section.

\subsection{Cost structure}\label{subsec:costs}
We introduce the three cost structures \textrm{C1}, \textrm{C2} and \textrm{C3}, which are presented in \mbox{Table \ref{tab:cost-ratio}}. To discourage repositioning tasks that yield negligible gain, we introduce a small travel cost which is paid each time unit an engineer is traveling. Each cost structure represents a distinct, realistic relationship between preventive and corrective costs that induces distinctive optimal policies favoring more or less frequent maintenance actions. For example, when $ \nicefrac{c^{\textrm{CM}} }{c^{\textrm{PM}} }$ is large, i.e., when CM costs greatly surpass PM costs, we expect that preventive maintenance policies outperform reactive policies. In all experiments, we consider a discount factor $\gamma = 0.99$.
\begin{table}[!ht]
\centering
\resizebox{0.55\textwidth}{!}{%
\begin{tabular}{c|cccc|c}
\toprule
Cost structure & $c^{\textrm{PM}}$ & $c^{\textrm{CM}}$ & $c^{\textrm{DT}}$ & $c^{\textrm{T}}$ & $\nicefrac{c^{\textrm{CM}}+c^{\textrm{DT}} }{c^{\textrm{PM}}+c^{\textrm{DT}} }$ \\ \midrule
\textrm{C1} & 0 & 0 & 1 & 0.05 & 1 \\
\textrm{C2} & 1 & 2 & 10 & 0 & 1.09 \\
\textrm{C3} & 1 & 4 & 1 & 0.05 & 2.5 \\
 \bottomrule
\end{tabular}}
\caption{Cost structures considered in the numerical experiments.}
\label{tab:cost-ratio}
\end{table}
\subsection{Hospital networks}\label{subsec: networks}
We construct six $K$-DTMPA instances, each having a different combination of network size, cost structure and machine degradation matrices. Besides the networks introduced in \cite{dtmpa}, we construct two additional geographical layouts with real-life asset network characteristics. 

To this end, these latter layouts are based on the Dutch hospital network. This case is appropriate since hospital equipment includes medical imaging and image-guided therapy systems. Such systems are associated with high costs, and manufacturers of such systems increasingly seek to avoid unplanned downtime via remote monitoring. The corresponding travel time matrices are constructed using the four-digit zip codes found in the 2021 hospital inventory data set published by the Dutch National Institute for Public Health and the Environment (RIVM). The partial zip codes are converted to GPS coordinates using the 4PP data set maintained by the Dutch Key Register Addresses and Buildings (Basisregistratie Adressen en Gebouwen). Finally, using these GPS coordinates, we compute the travel time between locations using the public OpenStreetMap application programming interface. The obtained travel times are converted and rounded up to multiples of 15 minutes. We assume repairs take $4$ time units, i.e., $1$ hour, regardless of the machine, the engineer or the maintenance type.

Generally, we consider two types of degradation matrices: One type with only two states, used to create the D\&R instances, and one type with three states, to include the aspect of preventive maintenance which increases the problem complexity significantly. We briefly discuss specifics regarding the instances. 

\subsubsection*{Academic hospitals}
\begin{figure}[!ht]
 \centering
 \begin{minipage}{.5\textwidth}
 \centering
 \includegraphics[width=0.8\linewidth]{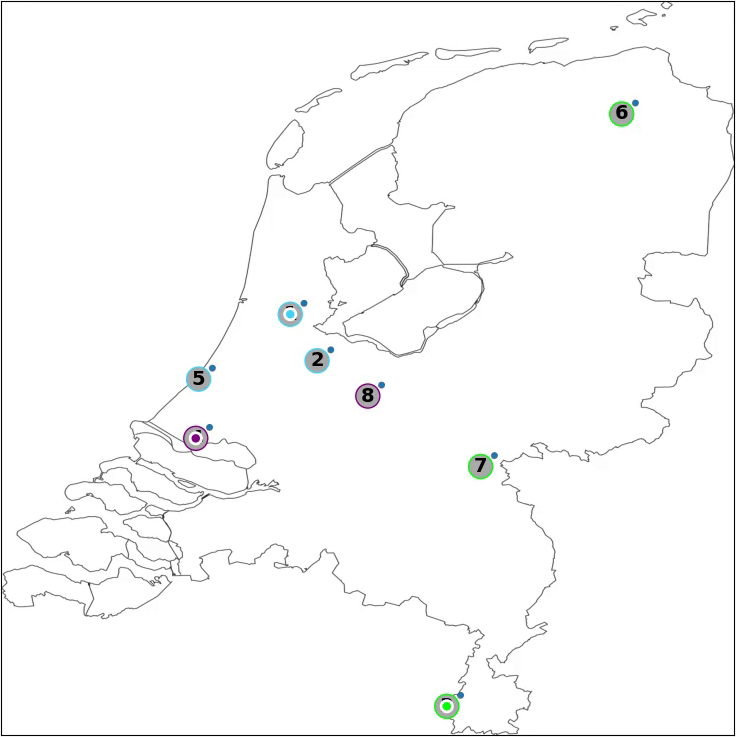}
 \end{minipage}%
 \begin{minipage}{0.5\textwidth}
 \centering
\resizebox{0.95\columnwidth}{!}{%
\begin{tabular}{c|llllllll}
$\theta$ & \rotatebox{90}{Amsterdam } & \rotatebox{90}{Amsterdam } & \rotatebox{90}{Maastricht} & \rotatebox{90}{Rotterdam} & \rotatebox{90}{Leiden} & \rotatebox{90}{Groningen} & \rotatebox{90}{Nijmegen} & \rotatebox{90}{Utrecht} \\
\hline
Amsterdam & 0 & 1 & 11 & 4 & 3 & 10 & 7 & 3 \\
Amsterdam & 1 & 0 & 11 & 5 & 3 & 10 & 7 & 3 \\
Maastricht & 11 & 11 & 0 & 11 & 12 & 17 & 8 & 10 \\
Rotterdam & 4 & 5 & 11 & 0 & 3 & 13 & 7 & 4 \\
Leiden & 3 & 3 & 12 & 3 & 0 & 12 & 8 & 4 \\
Groningen & 10 & 10 & 17 & 13 & 12 & 0 & 11 & 10 \\
Nijmegen & 7 & 7 & 8 & 7 & 8 & 11 & 0 & 5 \\
Utrecht & 3 & 3 & 10 & 4 & 4 & 10 & 5 & 0 \\
\end{tabular}
}
 \end{minipage}
 \caption{(Figure best viewed in color.) The Dutch academic hospitals with the corresponding travel time matrix $\Theta$ in quarters. The engineers are colored \textit{cyan}, \textit{green} and \textit{purple} and are located in Amsterdam, Maastricht and Rotterdam, respectively. For the decomposition heuristic, appropriate clusters are constructed using $K$-means clustering; locations within the respective clusters of engineers are colored accordingly.}
 \label{fig:academic_hospitals}
\end{figure}

The first two cases contain the $8$ academic hospitals in the Netherlands and will serve as complex yet relatively understandable and manageable examples to study the behavior of the learned policies. The Dutch academic hospitals are located in Amsterdam (2x), Groningen, Leiden, Maastricht, Nijmegen, Rotterdam and Utrecht, and are serviced by $K=3$ maintenance engineers, see \mbox{Figure \ref{fig:academic_hospitals}} for a visualization of the geographical layout and the corresponding travel time matrix. For the D\&R $K$-DTMPA instance, we adopt the degradation matrix \textrm{$\tilde{\textrm{Q}}$1} together with cost structure \textrm{C1}, referred to as M8K3-$\tilde{\textrm{Q}}$1C1. Under these dynamics, few failures occur and the objective is to minimize machine unavailability. \mbox{Figure \ref{fig:academic_hospitals}} shows a network decomposition of M8K3-$\tilde{\textrm{Q}}$1C1 into $1$-DTMPA instances. Each engineer induces a $1$-DTMPA instance on the locations within their assigned cluster, e.g., the first engineer maintains the locations in Amsterdam (2x) and Leiden (and only those).
\[\small{
\textrm{$\tilde{\textrm{Q}}$1} = \begin{bmatrix}
\frac{199}{200} & \frac{1}{200}\\
0 & 1 \\
\end{bmatrix}\quad
\textrm{$\tilde{\textrm{Q}}$2} = \begin{bmatrix}
\frac{149}{150} & \frac{1}{150} & 0\\
0 & \frac{49}{50} & \frac{1}{50} \\
0 & 0 & 1
\end{bmatrix}\quad
\label{mat:q2}}
\]
For the preventive maintenance $K$-DTMPA instance M8K3-$\tilde{\textrm{Q}}$2C3, we adopt the degradation matrix \textrm{$\tilde{\textrm{Q}}$2} together with cost structure \textrm{C3}. This models the situation where the alert is issued on average at $75\%$ of the machine's life expectancy and is thus an accurate indicator of failure. The goal now is to perform preventive maintenance while keeping all the machines operational.

\subsubsection*{City hospitals}
\begin{figure}[!ht]
 \centering
 \begin{minipage}{.5\textwidth}
 \centering
 \includegraphics[width=0.8\linewidth]{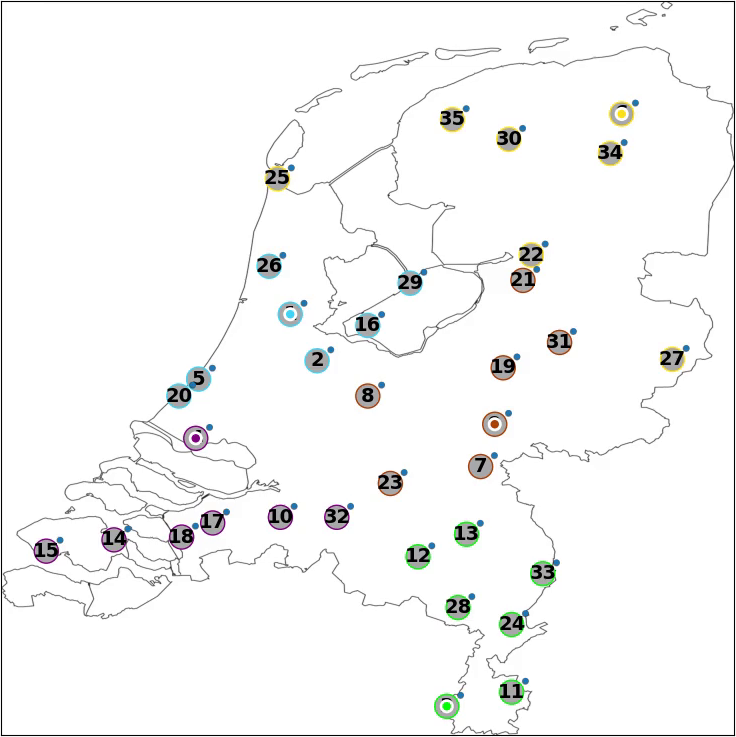}
 \end{minipage}%
 \begin{minipage}{0.5\textwidth}
 \centering
%\[\small{
\begin{align*}
&\textrm{$\tilde{\textrm{Q}}$3} = \begin{bmatrix}
\frac{399}{400} & \frac{1}{400}\\
0 & 1 \\
\end{bmatrix}\quad
\\
&\textrm{$\tilde{\textrm{Q}}$4} = \begin{bmatrix}
\frac{299}{300} & \frac{1}{300} & 0\\
0 & \frac{99}{100} & \frac{1}{100} \\
0 & 0 & 1
\end{bmatrix}\quad
\end{align*}
%\]
 \end{minipage}
 \caption{The subset of Dutch city hospitals with the corresponding degradation matrices \textrm{$\tilde{\textrm{Q}}$3} and \textrm{$\tilde{\textrm{Q}}$4}. The color-coding is as in \mbox{Figure \ref{fig:academic_hospitals}}, the additional engineers are colored \textit{brown} and \textit{yellow} and are located in Arnhem and Groningen, respectively. For the decomposition heuristic, appropriate clusters are constructed using $K$-means clustering; locations within the respective clusters of engineers are colored accordingly.}
 \label{fig:city_hospitals}
\end{figure}

We extend the academic hospital network by including a geographically dispersed subset of $35$ city hospitals. The network is serviced by $K=5$ maintenance engineers. This network will serve to study the behavior of the learned policies for industrial-scale $K$-DTMPA instances, see Figure \ref{fig:city_hospitals} for a visualization of the geographical layout and the asset degradation matrices. For the D\&R $K$-DTMPA instance, we adopt the degradation matrix \textrm{$\tilde{\textrm{Q}}$3} together with cost structure \textrm{C1}, referred to as M35K5-$\tilde{\textrm{Q}}$3C1. Under these dynamics, few failures occur and the objective is to minimize machine unavailability whilst maximizing coverage. For the preventive maintenance $K$-DTMPA instance M35K5-$\tilde{\textrm{Q}}$4C3, we adopt the degradation matrix \textrm{$\tilde{\textrm{Q}}$4} together with cost structure \textrm{C3}. Note that the alert is again issued (on average) at $75\%$ of the machine's life expectancy.

\section{Numerical results}\label{sec:numerical_results}
This section contains the experimental results for the previously introduced $K$-DTMPA instances. All values are obtained using $10^6$ repetitions, if applicable. The reported half-widths correspond to asymptotic $95\%$ confidence intervals. All DRL algorithm parameters used for training the neural networks are listed per instance in \ref{app: drl_param}. Information regarding the duration and the cost of training of the neural networks can be found in \ref{app: drl_training}. A complete list of experiments comparing the proposed heuristics acting on the same environment is available online\footnote{\url{https://retrospectiverotations.com/k-dtmpa/k-dtmpa.html}}.

\subsection{Impact of feature design on trained neural network policies}\label{subsec:feature_design}
First, we investigate the effectiveness of the proposed feature representation $f_3(h)$, cf. \mbox{Section \ref{subsec:feature_representation}}. Recall that $f_2(h)$ is obtained by removing the last entry of $f_1(h)$, i.e., $f_1(h) = \left( f_2(h), \sum_{m\in\mathcal{M}} n^\textrm{av}_m \right)$ and $f_3(h)$ denotes the most compact state representation. Given a data set $\mathcal{D}$ consisting of $500,000$ state-actions pairs corresponding to the reactive heuristic \mbox{$\pi_\textrm{D}$ $(s_m(h) \equiv x^\textbf{f}_m)$}, for each feature representation, we train $5$ neural network policies and report relevant statistics in \mbox{Table \ref{tab:feature_design}}.
\begin{table}[!ht]
\centering
\begin{tabular}{c|c|c|c|c}
\toprule
$f(h)$ & \textrm{min} $J(\pi_{\theta_1})$ & \textrm{max} $J(\pi_{\theta_1})$ & Average & CV \\ \midrule
$f_1(h)$ & $26.941 \pm 0.063$ & $27.013 \pm 0.063$ &$26.981$ & $1.011 \cdot 10^{-3}$ \\ \midrule
$f_2(h)$ & $26.953 \pm 0.063$ & $27.010 \pm 0.063$ &$26.974$ & $8.609 \cdot 10^{-4}$ \\ \midrule
$f_3(h)$ & $27.455 \pm 0.065$ & $28.024 \pm 0.067$ & $27.763$ & $9.017 \cdot 10^{-3}$\\ 
\bottomrule
\end{tabular}
\caption{One-step policy improvement results for the dispatching \& repositioning instance M8K3-$\tilde{\textrm{Q}}$1C1 highlighting the variability of the trained neural network when varying the feature design $f(h)$. In all cases, the initial policy is the reactive dispatching heuristic \mbox{$\pi_\textrm{D}$ $(s_m(h) \equiv x^\textbf{f}_m)$}. We train $5$ neural network policies per choice of $f(h)$ and report the performance of the best and worst policy, as well as the average and coefficient of variation (CV) of the acquired performance estimates. The benchmark satisfies $J(\pi_\textrm{D}\text{ }(s_m(h) \equiv x^\textbf{f}_m) ) = 27.612 \pm 0.065$.}
\label{tab:feature_design}
\end{table}
The choice of the feature design proves crucial: trained neural network policies using feature design $f_3(h)$ (the most compact state representation) barely beat the benchmark set by the reactive heuristic, if at all. Using the proposed feature design $f_1(h)$ results in less training variability and consistently produces policies that perform significantly better. This also holds for the feature representation $f_2(h)$, therefore, dropping the last feature does not significantly affect performance. Additional experiments for $f_1(h)$ where we also varied the data set $\mathcal{D}$ produced similar results. In the forthcoming sections, in all experiments, we use the feature representation $f_1(h)$.

\subsection{Single maintainer instances}\label{subsec: opt_gap}
To demonstrate that DCL can produce the optimal policy, we focus on the DTMPA instances M4-Q2Q3 and M6-Q2Q3Q4 introduced by \citet[Section 6.3]{dtmpa}, under cost structure C2 listed in \mbox{Table \ref{tab:cost-ratio}}. Note that $c^\textrm{T} = 0$, i.e., there is no cost for travel. Moreover, repair times and travel times are assumed to be $1$. The corresponding Q-matrices differ from the $\tilde{\textrm{Q}}$-matrices introduced in \mbox{Section \ref{subsec: networks}}, and can be found in \mbox{\ref{app: matrices}}. The shorthand notation must be interpreted as follows: First, the number of machines and engineers is listed, followed by a sequence of degradation matrices which are assumed to be distributed evenly over the machines. For example, the $1$-DTMPA instance M4K1-Q2Q3C2 contains four machines: Two with matrix Q2 and two with matrix Q3, all sharing cost structure C2.

\paragraph{M4K1-Q2Q3C2}
For small instances, the optimal policy $\pi_*$ can be obtained via exact policy iteration. We perform two policy improvement steps using DCL on three dispatching heuristics $\pi_D$, with maintenance thresholds $s_m(h)$ ranging from $3$ to $5$, and on both the random policy $\pi_\textrm{R}$ and the idle policy $\pi_\textrm{I}$. DCL consistently produces a near-optimal policy after only two iterations, regardless of the initial solution. The best-found neural network policy places the optimality gap (computed as the relative increase over $J(\pi_*)$) at only $0.45\%$. The improvement over the $n$QR-DDQN policy proposed by \citet{dtmpa} is $8.34\%$. (Note that DCL operates at an information advantage compared to $n$QR-DDQN, since the latter only observes transitions to the first degradation state $x^{\textbf{h}}+1$ and the failed state $x^{\textbf{f}}$, while the former observes all degradation state transitions.)

Observe from \mbox{Table \ref{tab:result-M4K1-Q2Q3C2}} that the best performing dispatching heuristic ($\pi_D$ with $s_m(h) \equiv 3$) is not necessarily the best choice of initial policy for a single one-step improvement; the neural network policy trained using the reactive dispatching heuristic (with $s_m(h) \equiv x^\textbf{f}_m \equiv 5$) performs at least $1.58\%$ better than any of the other first-generation neural network policies. The idle policy is a particularly poor choice since it produces a rather homogeneous data set and only learns to start maintenance at the current location in the first iteration. 
\begin{table}[!ht]
\centering
\resizebox{1\linewidth}{!}{%
\begin{tabular}{c|c|c|c|c|c}
\toprule
$\pi$ & $\pi_\textrm{D}$ $(s_m(h) \equiv 3)$ & $\pi_\textrm{D}$ $(s_m(h) \equiv 4)$ & $\pi_\textrm{D}$ $(s_m(h) \equiv x^\textbf{f}_m)$ & $\pi_\textrm{R}$ & $\pi_\textrm{I}$ \\ \midrule
$J(\pi)$ & $659.914 \pm 1.380$ & $599.654 \pm 1.243$ & $ 780.818 \pm 1.631$ & $2313.600 \pm 5.006 $ & $3509.960 \pm 7.732$\\ \midrule
$J(\pi_{\theta_{1}})$ & $490.577 \pm 1.000$ & $453.732 \pm 0.931$ & $446.682 \pm 0.919$ & $515.649 \pm 1.039 $ & $2833.530 \pm 6.156 $ \\ 
\midrule
$J(\pi_{\theta_{2}})$ & $435.047 \pm 0.896$ & $434.385 \pm 0.896$ & $434.393 \pm 0.896$ & $440.285 \pm 0.909$ & $1229.250 \pm 2.634$ \\ 
\bottomrule
\end{tabular}
}
\caption{One-step policy improvement results for the single maintainer instance M4K1-Q2Q3C2. The optimal solution satisfies $J(\pi_*) = 432.440$ \citep[Table 2]{dtmpa}.}
\label{tab:result-M4K1-Q2Q3C2}
\end{table}

\paragraph{M6K1-Q2Q3Q4C2}
For this instance, computing the optimal policy is intractable and thus the best available solution in literature is the neural network policy trained by $n$QR-DDQN. We perform three policy improvement steps using DCL on three dispatching heuristics $\pi_D$ with varying maintenance thresholds and the random policy $\pi_\textrm{R}$. Similarly, the best neural network policy yields a $12.34\%$ advantage over the $\pi_{\textrm{$n$QR-DDQN}}$ policy (which may be partially due to informational advantage) and the reactive dispatching heuristic consistently produces the best policy after the first iteration. The results on the M4K1-Q2Q3C2 instance, together with observations from \mbox{Table \ref{tab:result-M6K1-Q2Q3Q4C2}}, indicate that the neural network policy improvements have ended which could be taken as weak evidence that the best performing neural network policy is near-optimal.
\begin{table}[!ht]
\centering
\begin{tabular}{c|c|c|c|c}
\toprule
$\pi$ & $\pi_\textrm{D}$ $(s_m(h) \equiv 4)$ & $\pi_\textrm{D}$ $(s_m(h) \equiv 5)$ & $\pi_\textrm{D}$ $(s_m(h) \equiv x^\textbf{f}_m)$ & $\pi_\textrm{R}$ \\ \midrule
$J(\pi)$ & $1100.490 \pm 2.368 $ & $1129.070 \pm 2.391$ & $1207.200 \pm 2.572$ & $3876.330 \pm 8.644$ \\ \midrule
$J(\pi_{\theta_{1}})$ & $683.770 \pm 1.423$ & $704.293 \pm 1.494$ & $658.454 \pm 1.386$ & $ 765.695 \pm 1.588$ \\ 
\midrule
$J(\pi_{\theta_{2}})$ & $625.079 \pm 1.309$ & $635.208 \pm 1.329$ & $628.127 \pm 1.317$ & $641.696 \pm 1.343$ \\ \midrule 
$J(\pi_{\theta_{3}})$ & $ 623.407 \pm 1.305$ & $ 623.793 \pm 1.307$ & $ 623.593 \pm 1.305$ & $623.596 \pm 1.305$ \\ 
\bottomrule
\end{tabular}
\caption{One-step policy improvement results for the single maintainer instance M6K1-Q2Q3Q4C2. The estimate for the neural network policy $\pi_{\textrm{$n$QR-DDQN}}$ satisfies $J(\pi_{\textrm{$n$QR-DDQN}}) = 711.188 \pm 5.398$ \citep[Table 2]{dtmpa}.}
\label{tab:result-M6K1-Q2Q3Q4C2}
\end{table}
\subsection{Dispatching \& repositioning instances}\label{subsec: one_step_API}
To illustrate how DCL improves upon an existing solution, we now turn our attention to the $K$-DTMPA instances M8K3-$\tilde{\textrm{Q}}$1C1 and M35K5-$\tilde{\textrm{Q}}$3C1. For both instances, the benchmark is set by the reactive dispatching heuristic \mbox{$\pi_\textrm{D}$ $(s_m(h) \equiv x^\textbf{f}_m)$}. We note that for D\&R instances, this benchmark is expected to be rather strong and difficult to beat.

\paragraph{M8K3-$\tilde{\textrm{Q}}$1C1}
In this $3$-DTMPA instance, the engineers are initially placed in the cities Amsterdam, Maastricht and Rotterdam. We perform three policy improvement steps using DCL on the reactive dispatching heuristic \mbox{$\pi_\textrm{D}$ $(s_m(h) \equiv x^\textbf{f}_m)$}, the random policy $\pi_\textrm{R}$ and the decomposition heuristic $\pi^\textrm{DEC}_\textrm{D}$, for which we show the results in \mbox{Table \ref{tab:result-M8K3-Q1C1}}. The best neural network policy $\pi_{\theta_1}$ improves upon the benchmark with $1.74\%$. A second iteration yields a further $0.67\%$ performance improvement and a third step accomplishes an additional $1.07\%$ cost reduction and thus a total performance improvement of $3.48\%$ over the benchmark. Three steps of DCL improving $\pi_\textrm{R}$ is not sufficient to outperform the benchmark, meaning that $\pi_\textrm{R}$ is not a suitable initial solution in a cooperative setting. Moreover, decomposing the network into clusters and solving the induced $1$-DTMPA instances individually does not bring us close to the benchmark either. 

\begin{table}[!ht]
\centering
\begin{tabular}{c|c|c|c}
\toprule
$\pi$ & $\pi_\textrm{D}$ $(s_m(h) \equiv x^\textbf{f}_m)$ & $\pi_\textrm{R}$ & $\pi^\textrm{DEC}_\textrm{D}$ $(s_m(h) \equiv x^\textbf{f}_m)$ \\ \midrule
$J(\pi)$ & $27.612 \pm 0.065$ & $218.390 \pm 0.649$ & $33.143 \pm 0.067$ \\ \midrule
$J(\pi_{\theta_1})$ & $27.131 \pm 0.067$ & $64.605 \pm 0.143$ & $30.779 \pm 0.055 $ \\ \midrule
$J(\pi_{\theta_2})$ & $26.945 \pm 0.065$ & $51.320 \pm 0.112$ & $30.727 \pm 0.055
$\\ \midrule
$J(\pi_{\theta_3})$ & $26.652 \pm 0.061$ & $46.774 \pm 0.102$ & $30.727 \pm 0.055$\\ \bottomrule
\end{tabular}
\caption{One-step policy improvement results for the dispatching \& repositioning instance M8K3-$\tilde{\textrm{Q}}$1C1. Note that the results for the reactive dispatching heuristic are obtained using a data set $\mathcal{D}$ that is smaller than the data set used for the experiments in \mbox{Table \ref{tab:feature_design}}.}
\label{tab:result-M8K3-Q1C1}
\end{table}
\begin{figure}[!ht]
 \centering
 \begin{subfigure}[t]{0.32\textwidth}
 \centering
 \includegraphics[width=\textwidth]{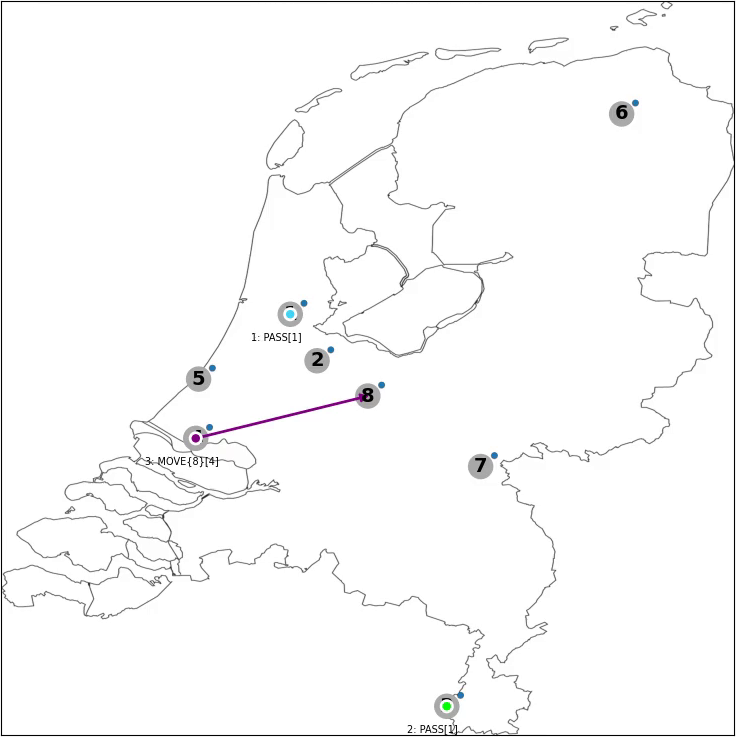}
 \caption{Strategic initial positioning.}
 \label{subfig:M8K3-Q1C1_1}
 \end{subfigure}%
 \hfill
 \begin{subfigure}[t]{0.32\textwidth}
 \centering
 \includegraphics[width=\textwidth]{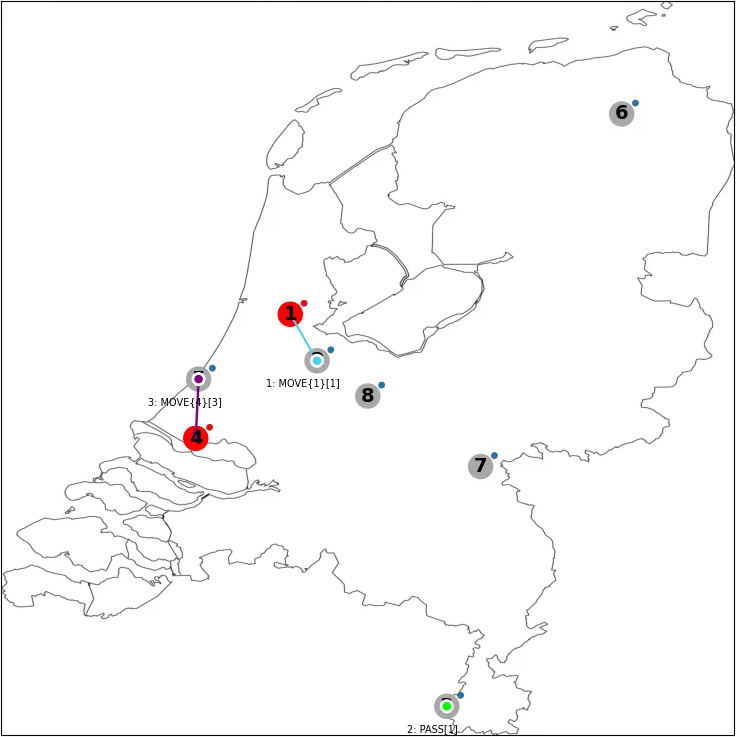}
 \caption{Efficient dispatching.}
 \label{subfig:M8K3-Q1C1_2}
 \end{subfigure}
 \hfill
 \begin{subfigure}[t]{0.32\textwidth}
 \centering
 \includegraphics[width=\textwidth]{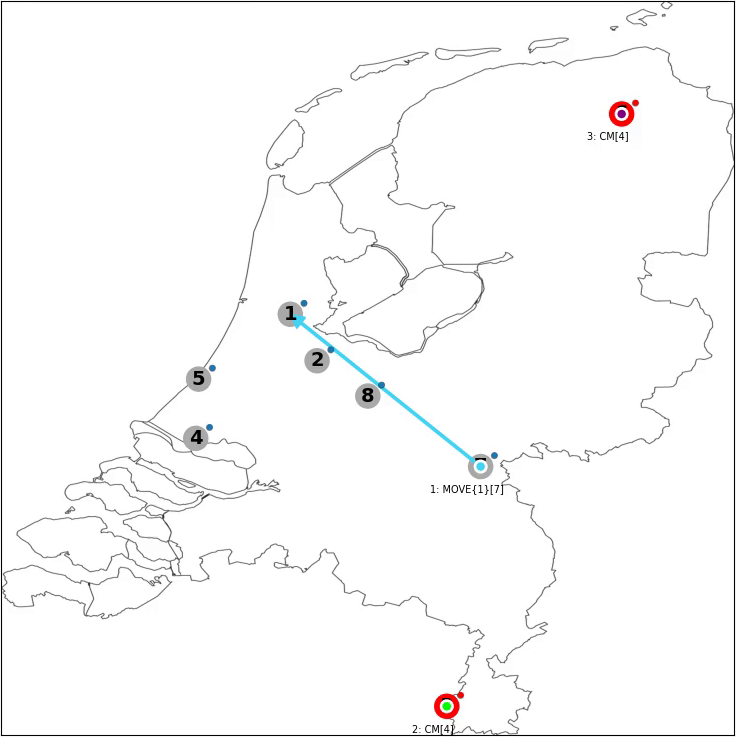}
 \caption{Tactical repositioning.}
 \label{subfig:M8K3-Q1C1_3}
 \end{subfigure}
\caption{(Figure best viewed in color.) Policy aspects of the DCL improved reactive dispatching heuristic for the dispatching \& repositioning instance M8K3-$\tilde{\textrm{Q}}$1C1. The color-coding is as in \mbox{Figure \ref{fig:city_hospitals}}. The labels \textsc{pass}, \textsc{move}, \textsc{pm} and \textsc{cm}, correspond to wait, move to another location and preventive/corrective maintenance actions. \textit{Blue} dots on top of machine nodes indicate that the machine is healthy, \textit{orange} when alerted or \textit{red} when the machine is down.}
 \label{fig:api_policy}
\end{figure}

We next briefly illustrate how the best neural network policy improves upon the benchmark in terms of the behavioral aspects introduced in \mbox{Section \ref{subsec:policy_aspects}}. In \mbox{Figure \ref{subfig:M8K3-Q1C1_1}}, we observe that the DRL agent moves an available engineer from Rotterdam to the centrally located Utrecht, which is likely a better initial placement of the engineer. The DRL agent's dispatching strategy is similar to the reactive heuristic and handles the tricky cases correctly as well, see for instance the dispatching problem in \mbox{Figure \ref{subfig:M8K3-Q1C1_2}}. When there are events at remote locations in the network, the DRL agent has learned to proactively move an engineer to achieve better coverage. \mbox{Figure \ref{subfig:M8K3-Q1C1_3}} shows such a tactical repositioning of an engineer from Nijmegen to Amsterdam.

\begin{figure}[!ht]
\begin{subfigure}{.5\textwidth}
 \centering
 % include first image
 \includegraphics[width=0.8\linewidth]{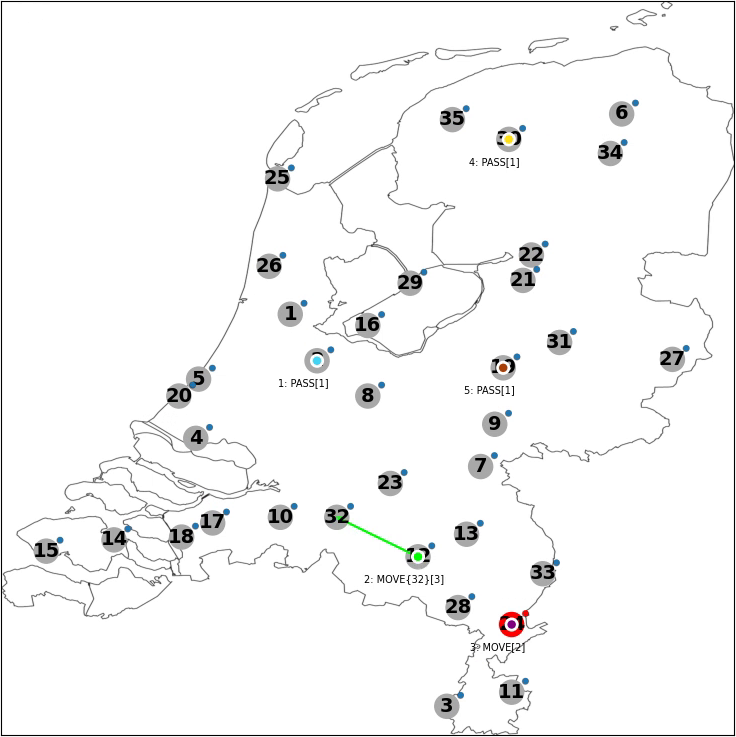} 
 \caption{Tactical repositioning from Eindhoven to Tilburg.}
 \label{subfig:M35K5-Q3C1_1}
\end{subfigure}
\begin{subfigure}{.5\textwidth}
 \centering
 % include second image
 \includegraphics[width=0.8\linewidth]{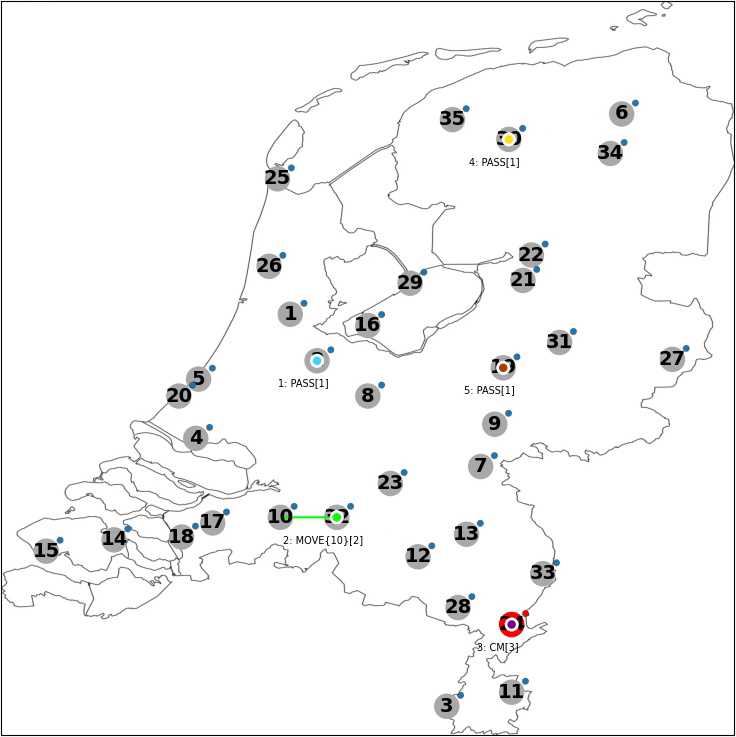} 
 \caption{Tactical repositioning from Tilburg to Breda.}
 \label{subfig:M35K5-Q3C1_2}
\end{subfigure}
\caption{Policy aspects of the DCL improved reactive dispatching heuristic for the dispatching \& repositioning instance M35K5-$\tilde{\textrm{Q}}$3C1. The color-coding is as in \mbox{Figure \ref{fig:api_policy}}.}
\label{fig:api_policy_2}
\end{figure}

\paragraph{M35K5-$\tilde{\textrm{Q}}$3C1}
In the first large-scale instance, the initial placement of the additional two engineers is Arnhem and Groningen. We perform two policy improvement steps using DCL on the reactive dispatching heuristic \mbox{$\pi_\textrm{D}$ $(s_m(h) \equiv x^\textbf{f}_m)$} and the decomposition heuristic $\pi^\textrm{DEC}_\textrm{D}$, for which we show the results in \mbox{Table \ref{tab:result-M35K5-Q3C1}}. 
\begin{wraptable}{r}{0.5\linewidth}
\begin{tabular}{c|c|c}
\toprule
$\pi$ & $\pi_\textrm{D}$ $(s_m(h) \equiv x^\textbf{f}_m)$ & $\pi^\textrm{DEC}_\textrm{D}$ $(s_m(h) \equiv x^\textbf{f}_m)$ \\ \midrule
$J(\pi)$ & $65.612 \pm 0.145$ & $71.590 \pm 0.090$\\ \midrule
$J(\pi_{\theta_1})$ & $63.178 \pm 0.137$ & $67.134 \pm 0.120$ \\ \midrule
$J(\pi_{\theta_2})$ & $62.101 \pm 0.135$ & $66.326 \pm 0.104$ \\ \bottomrule
\end{tabular}
\caption{One-step policy improvement results for the \\ dispatching \& repositioning instance M35K5-$\tilde{\textrm{Q}}$3C1.}
\label{tab:result-M35K5-Q3C1}
\end{wraptable}
The best found neural network policy $\pi_{\theta_2}$ improves upon the benchmark with $5.35\%$. Besides a better strategic initial positioning, the neural network policy learns to travel using intermediate locations, see \mbox{Figure \ref{fig:api_policy_2}}. The benefit of such behavior is that it enables the engineer to have an extra decision epoch, e.g., to divert to another location or even catch a failure at one of the intermediate locations. Moreover, after two iterations of DCL, sharing resources over the network is shown to yield a $6.37\%$ cost improvement.

\subsection{Preventive maintenance instances}\label{subsec:preventive_maintenance_instances}
 We show that DCL can also handle the more complex $K$-DTMPA instances M8K3-$\tilde{\textrm{Q}}$2C3 and M35K5-$\tilde{\textrm{Q}}$4C3. For these instances, no strong benchmark exists and therefore, we consider both the greedy and the reactive dispatching heuristics $\pi_\textrm{D}$, as well as the decomposition heuristic $\pi^\textrm{DEC}_\textrm{D}$.

\paragraph{M8K3-$\tilde{\textrm{Q}}$2C3}
We perform three steps of DCL on the various heuristics, for which we show the results in \mbox{Table \ref{tab:result-M8K3-Q2C3}}. The best neural network policy $\pi_{\theta_3}$ improves upon the best heuristic with $6.44\%$ and shows that at least a $5.82\%$ cost reduction can be achieved by sharing resources over the network. We briefly illustrate how the third generation neural network policy $\pi_{\theta_3}$ improves upon the benchmark in terms of the behavioral aspects introduced in \mbox{Section \ref{subsec:policy_aspects}}.
\begin{table}[!ht]
\centering
\begin{tabular}{c|c|c|c|c}
\toprule
$\pi$ & $\pi_\textrm{D}$ $(s_m(h) \equiv 2)$ & $\pi_\textrm{D}$ $(s_m(h) \equiv x^\textbf{f}_m)$ & $\pi^\textrm{DEC}_\textrm{D}$ $(s_m(h) \equiv 2)$ & $\pi^\textrm{DEC}_\textrm{D}$ $(s_m(h) \equiv x^\textbf{f}_m)$ \\ \midrule
$J(\pi)$ & $26.736 \pm 0.061$ & $31.756 \pm 0.090$ & $29.154 \pm 0.051$ & $35.538 \pm 0.074$\\ \midrule
$J(\pi_{\theta_{1}})$ & $25.608 \pm 0.065$ & $26.179 \pm 0.069$ & $27.141 \pm 0.053$ & $27.159 \pm 0.047$ \\ \midrule
$J(\pi_{\theta_{2}})$ & $25.275 \pm 0.067$ & $25.455 \pm 0.069$ & $26.805 \pm 0.051$ & $26.571 \pm 0.049$\\ \midrule
$J(\pi_{\theta_{3}})$ & $25.068 \pm 0.067$ & $25.148 \pm 0.067$ & $26.517 \pm 0.047$ & $26.571 \pm 0.049$ \\ 
\bottomrule
\end{tabular}
\caption{One-step policy improvement results for the preventive maintenance instance M8K3-$\tilde{\textrm{Q}}$2C3.}
\label{tab:result-M8K3-Q2C3}
\end{table}
\begin{figure}[!ht]
 \centering
 \begin{subfigure}[t]{0.32\textwidth}
 \centering
 \includegraphics[width=\textwidth]{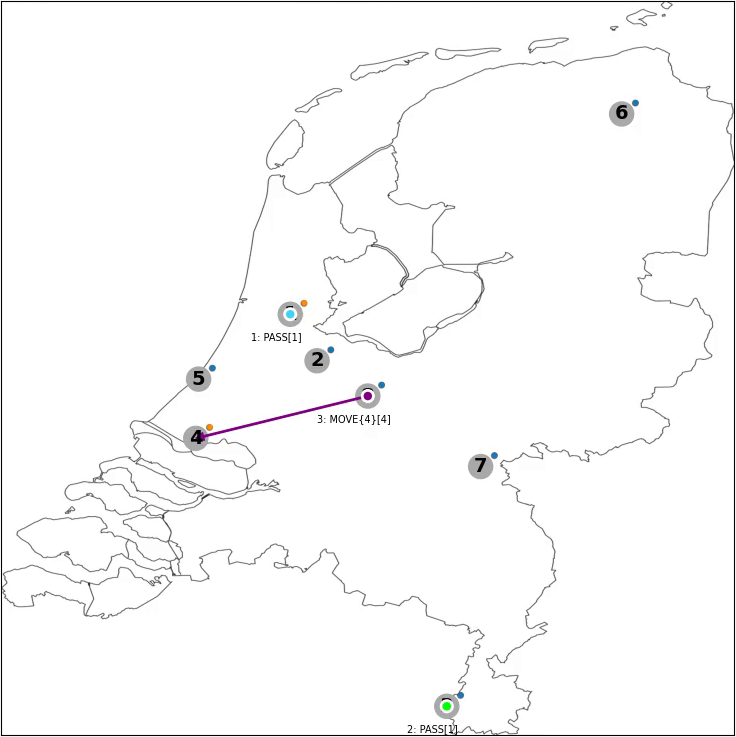}
 \caption{The trained agent proactively dispatches engineers to alerted locations.}
 \label{subfig:M8K3-Q2C3_1}
 \end{subfigure}%
 \hfill
 \begin{subfigure}[t]{0.32\textwidth}
 \centering
 \includegraphics[width=\textwidth]{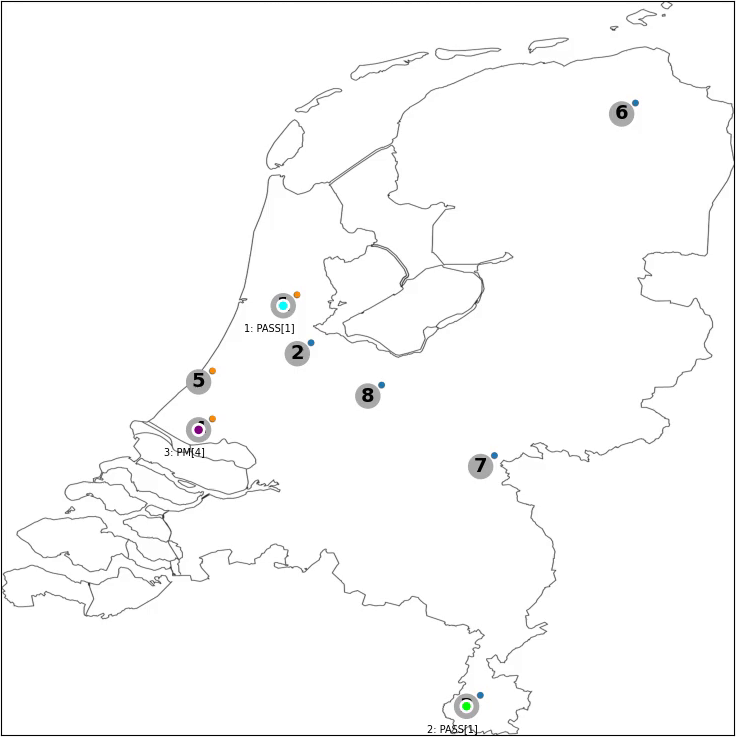}
 \caption{The trained agent performs preventive maintenance when there are many events in the network.}
 \label{subfig:M8K3-Q2C3_2}
 \end{subfigure}
 \hfill
 \begin{subfigure}[t]{0.32\textwidth}
 \centering
 \includegraphics[width=\textwidth]{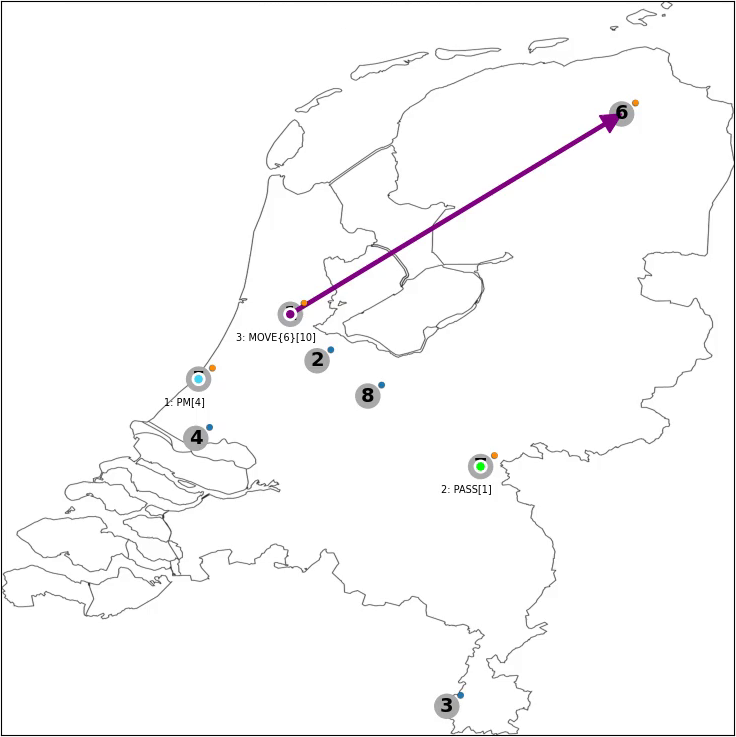}
 \caption{The trained agent tactically \\repositions engineers when there are \\distant events in the network.}
 \label{subfig:M8K3-Q2C3_3}
 \end{subfigure}
\caption{Policy aspects of the DCL improved greedy dispatching heuristic for the preventive maintenance instance M8K3-$\tilde{\textrm{Q}}$2C3. The color-coding is as in \mbox{Figure \ref{fig:api_policy}}.}
 \label{fig:api_policy_3}
\end{figure}
In Figure \ref{subfig:M8K3-Q2C3_1}, we see that the trained agent proactively moves available engineers to alerted locations, but postpones maintenance on them. When there are many alerts in the network, the trained agent chooses to initiate preventive maintenance, see \mbox{Figure \ref{subfig:M8K3-Q2C3_2}}. When there are events at remote locations in the network as \mbox{Figure \ref{subfig:M8K3-Q2C3_3}}, the DRL agent proactively dispatches an engineer.

\begin{figure}[!ht]
\begin{subfigure}{.5\textwidth}
 \centering
 % include first image
 \includegraphics[width=0.8\linewidth]{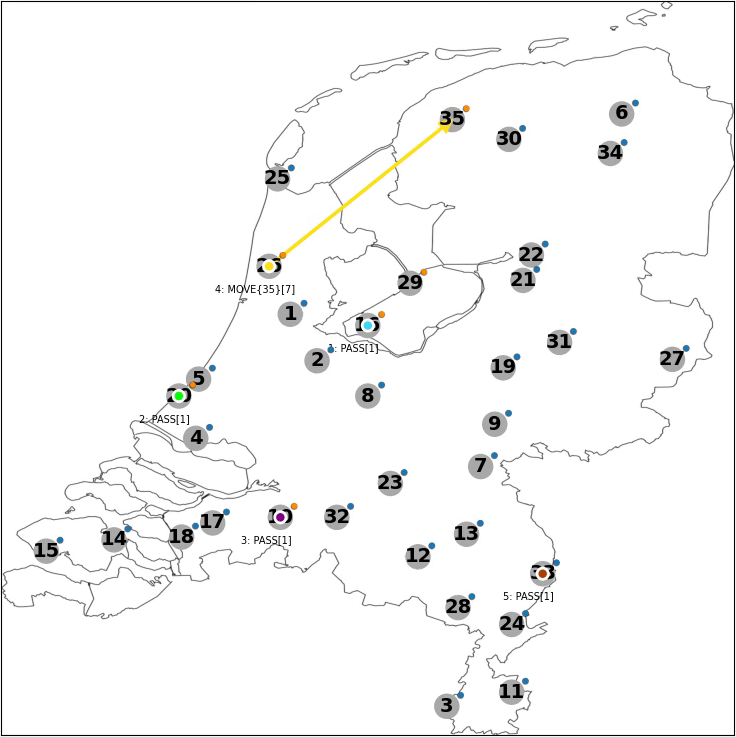} 
 \caption{The trained agent proactively moves engineers to \\alerted locations, ensuring a large coverage of the network.}
 \label{subfig:M35K5-Q4C3_1}
\end{subfigure}
\begin{subfigure}{.5\textwidth}
 \centering
 % include second image
 \includegraphics[width=0.8\linewidth]{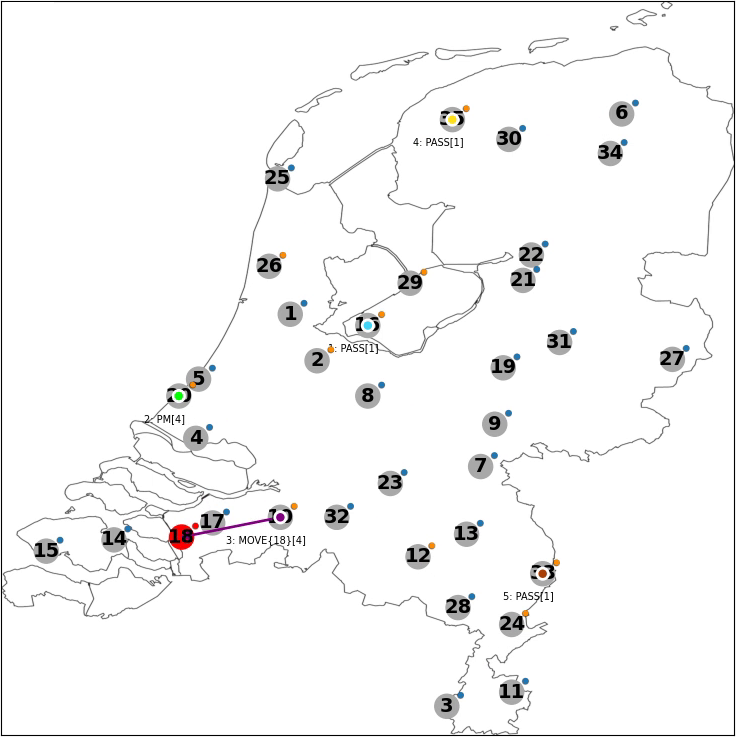} 
 \caption{The trained agent prioritizes moving to and \\ starting maintenance on assets in the failed state.}
 \label{subfig:M35K5-Q4C3_2}
\end{subfigure}
\caption{Policy aspects of the DCL improved reactive dispatching heuristic for the preventive maintenance instance M35K5-$\tilde{\textrm{Q}}$4C3. The color-coding is as in \mbox{Figure \ref{fig:api_policy}}.}
\label{fig:api_policy_4}
\end{figure}

\paragraph{M35K5-$\tilde{\textrm{Q}}$4C3}
We perform a single iteration of DCL on the reactive dispatching heuristic \mbox{$\pi_\textrm{D}$ $(s_m(h) \equiv x^\textbf{f}_m)$} due to the relatively large cost of additional steps (see \ref{app: drl_training}). We compare against three policy improvement steps on the network decomposition heuristics $\pi^\textrm{DEC}_\textrm{D}$, for which we show the results in \mbox{Table \ref{tab:result-M35K5-Q4C3}}. The neural network policy $\pi_{\theta_1}$ improves upon the initial policy with $6.08\%$ and shows that a $1.17\%$ cost reduction can be achieved by sharing resources over the network. In Figure \ref{subfig:M35K5-Q4C3_1}, we see that the trained agent proactively moves available engineers to alerted locations while ensuring a large coverage of the network. The trained agent however chooses to prioritize corrective maintenance, see \mbox{Figure \ref{subfig:M35K5-Q4C3_2}}. 
\begin{table}[!ht]
\centering
\begin{tabular}{c|c|c|c|c}
\toprule
$\pi$ & $\pi_\textrm{D}$ $(s_m(h) \equiv 2)$ & $\pi_\textrm{D}$ $(s_m(h) \equiv x^\textbf{f}_m)$ & $\pi^\textrm{DEC}_\textrm{D}$ $(s_m(h) \equiv 2)$ & $\pi^\textrm{DEC}_\textrm{D}$ $(s_m(h) \equiv x^\textbf{f}_m)$\\ \midrule
$J(\pi)$ & $57.512 \pm 0.125$ & $56.654 \pm 0.176$ & $58.365 \pm 0.067$ & $59.460 \pm 0.096$ \\ \midrule
$J(\pi_{\theta_{1}})$ & $-$ & $53.207 \pm 0.163$ & $58.262 \pm 0.092$ & $57.203 \pm 0.147 $ \\ \midrule
$J(\pi_{\theta_{2}})$ & $-$ & $-$ & $55.529 \pm 0.100 $ & $55.045 \pm 0.100 $ \\ \midrule
$J(\pi_{\theta_{3}})$ & $-$ & $-$ & $54.325 \pm 0.104 $ & $53.836 \pm 0.102 $ \\ \bottomrule
\end{tabular}
\caption{One-step policy improvement results for the preventive maintenance instance M35K5-$\tilde{\textrm{Q}}$4C3.}
\label{tab:result-M35K5-Q4C3}
\end{table}

\subsection{Robustness of policies to changes in the model}
In some cases, a neural network policy trained to optimize a $K$-DTMPA instance also yields a policy for a compatible $K$-DTMPA instance, but since the neural network was trained for a specific instance, this is somewhat detrimental to performance. We briefly investigate this next.

A trained neural network policy is also a policy for a compatible instance when the instance preserves the dimension of the feature vector $f(h)$, viz. when the number of machines remains the same or is reduced. (The case of removing machines can be tackled by replacing them with dummy machines that never emit alerts.) We investigate the effect of removing one hospital/machine (at location 2) and hiring/firing a single maintenance engineer (at location 2 and 1, respectively) for the cases M8K3-$\tilde{\textrm{Q}}$1C1, M35K5-$\tilde{\textrm{Q}}$3C1 and M8K3-$\tilde{\textrm{Q}}$2C3, for which we show the results in \mbox{Table \ref{tab:sensitivity_analysis}}. 
\begin{table}[!ht]
\tabcolsep=0.1cm
\resizebox{1\linewidth}{!}{%
\begin{tabular}{|c|cc|cc|cc|}
\hline
 $\mathbf{M-1}$ & \multicolumn{2}{c|}{M7K3-$\tilde{\textrm{Q}}$1C1} & \multicolumn{2}{c|}{M7K3-$\tilde{\textrm{Q}}$2C3} & \multicolumn{2}{c|}{M34K5-$\tilde{\textrm{Q}}$3C1} \\ \hline
$\pi$ & \multicolumn{1}{c|}{$\pi_\textrm{D}$ $(s_m(h) \equiv x^\textbf{f}_m)$} & $\pi_{\theta_3}$ $(s_m(h) \equiv x^\textbf{f}_m)$ & \multicolumn{1}{c|}{$\pi_\textrm{D}$ $(s_m(h) \equiv 2)$} & $\pi_{\theta_3}$ $(s_m(h) \equiv 2)$ & \multicolumn{1}{c|}{$\pi_\textrm{D}$ $(s_m(h) \equiv x^\textbf{f}_m)$} & $\pi_{\theta_2}$ $(s_m(h) \equiv x^\textbf{f}_m)$ \\ \hline 
$J(\pi)$ & \multicolumn{1}{c|}{$24.117 \pm 0.057$} & $23.908 \pm 0.055$ & \multicolumn{1}{c|}{ $23.390 \pm 0.053$ } & $21.833 \pm 0.061$ & \multicolumn{1}{c|}{$63.869 \pm 0.141$} & $61.085 \pm 0.133$ \\ \hline
 $\mathbf{K-1}$ & \multicolumn{2}{c|}{M8K2-$\tilde{\textrm{Q}}$1C1} & \multicolumn{2}{c|}{M8K2-$\tilde{\textrm{Q}}$2C3} & \multicolumn{2}{c|}{M35K4-$\tilde{\textrm{Q}}$3C1} \\ \hline
$\pi$ & \multicolumn{1}{c|}{$\pi_\textrm{D}$ $(s_m(h) \equiv x^\textbf{f}_m)$} & $\pi_{\theta_3}$ $(s_m(h) \equiv x^\textbf{f}_m)$ & \multicolumn{1}{c|}{$\pi_\textrm{D}$ $(s_m(h) \equiv 2)$} & $\pi_{\theta_3}$ $(s_m(h) \equiv 2)$ & \multicolumn{1}{c|}{$\pi_\textrm{D}$ $(s_m(h) \equiv x^\textbf{f}_m)$} & $\pi_{\theta_2}$ $(s_m(h) \equiv x^\textbf{f}_m)$ \\ \hline
 $J(\pi)$ & \multicolumn{1}{c|}{ $35.601 \pm 0.086$} & $98.050 \pm 0.417$ & \multicolumn{1}{c|}{$28.550 \pm 0.069$} & $29.584 \pm 0.110$ & \multicolumn{1}{c|}{$72.570 \pm 0.163$} & $69.442 \pm 0.153$ \\ \hline
$\mathbf{K+1}$ & \multicolumn{2}{c|}{M8K4-$\tilde{\textrm{Q}}$1C1} & \multicolumn{2}{c|}{M8K4-$\tilde{\textrm{Q}}$2C3} & \multicolumn{2}{c|}{M35K6-$\tilde{\textrm{Q}}$3C1} \\ \hline
$\pi$ & \multicolumn{1}{c|}{$\pi_\textrm{D}$ $(s_m(h) \equiv x^\textbf{f}_m)$} & $\pi_{\theta_3}$ $(s_m(h) \equiv x^\textbf{f}_m)$ & \multicolumn{1}{c|}{$\pi_\textrm{D}$ $(s_m(h) \equiv 2)$} & $\pi_{\theta_3}$ $(s_m(h) \equiv 2)$ & \multicolumn{1}{c|}{$\pi_\textrm{D}$ $(s_m(h) \equiv x^\textbf{f}_m)$} & $\pi_{\theta_2}$ $(s_m(h) \equiv x^\textbf{f}_m)$ \\ \hline
 $J(\pi)$ & \multicolumn{1}{c|}{ $23.477 \pm 0.055$ } & $118.617 \pm 0.314$ & \multicolumn{1}{c|}{ $26.202 \pm 0.059$ } & $54.196 \pm 0.163$ & \multicolumn{1}{c|}{$62.095 \pm 0.135$} & $64.428 \pm 0.141$ \\ \hline
\end{tabular}
}
\caption{Sensitivity analysis for the $K$-DTMPA instances M8K3-$\tilde{\textrm{Q}}$1C1, M35K5-$\tilde{\textrm{Q}}$3C1 and M8K3-$\tilde{\textrm{Q}}$2C3.}
\label{tab:sensitivity_analysis}
\end{table}
The case where only one machine is removed preserves the most similarity in the encountered features, the main difference being that one asset never emits alerts. In a network with identical assets, the trained neural network policy is thus robust against removing a single asset. Note that robustness against completely removing a machine is also an indication that performance would be robust against slight modification of machine degradation models. When removing a maintenance engineer, the trained neural network policy produced a good solution $2$ out of $3$ times. However, performance decreases significantly when adding a maintenance engineer, although the neural network policy still provides a suitable initial solution. An alternative approach to tackle such instances could be to construct a policy as follows: Use the trained neural network to dispatch $K$ engineers and employ a classical heuristic for the additional engineer. This, however, is outside the scope of this work.
\section{Conclusion and discussion}\label{sec:conclusion}
In this work, we study the dynamic traveling multi-maintainer problem with alerts ($K$-DTMPA) for a network of modern industrial assets with stochastic failure times maintained by $K$ maintenance engineers. We extended the existing $1$-DTMPA framework under the assumption of perfect information proposed by \cite{dtmpa}. Also, our experiments include cases with an underlying geographical nature, which are more challenging than the unit-distance cases considered by \cite{dtmpa}. In the $K$-DTMPA framework, independent degradation processes are observed in real-time by a central decision-maker. The decision-maker has access to perfect degradation information to decide on joint cost-effective dispatching, maintenance and repositioning actions for all available engineers.

To solve the problem, we adopt a deep reinforcement learning (DRL) approach based on approximate policy iteration, more specifically, deep controlled learning (DCL). To successfully apply DCL to $K$-DTMPA instances, we propose several new ideas: Actions for the engineers are selected sequentially and the feature design is tailored to each individual engineer. Moreover, we use policies tailored to the problem to kickstart DCL. This enables us to use expert knowledge without the requirement to penalize deviations from the expert-crafted policy (like \citet{DEMOOR2022535}). 

To demonstrate the effectiveness of our approach, we extended existing ranking heuristics to the multi-maintainer setting. More specifically, we equip ranking heuristics that rank alerts based on their observed degradation levels with a state-of-the-art dispatching algorithm. Moreover, we propose an additional benchmark heuristic through decomposition: We decompose the network into $K$ disjoint $1$-DTMPA instances using a handcrafted network clustering and solve each of the induced subproblems individually using DRL.

The results for small instances show that we can get close to the performance of optimal policies within a few iterations of the algorithm, regardless of the choice of the initial solution. When the problem complexity increases, the proposed DRL method yields effective policies that directly improve upon the benchmark. This significantly reduces the amount of required iterations, thereby saving costs. Moreover, by comparing with a traditional solution by decomposition, we show that it is cost-effective to share resources over the network.

Future research directions include expanding the proposed framework to include logistical and asset-maintainer constraints, possibly as a learning objective. The assumption of geometrically distributed degradation transition times remains and underlies the tractability of some of the instances; this assumption can be relaxed within the $K$-DTMPA framework, which we leave for future work. The bottleneck of the adopted DRL approach to solving large-scale $K$-DTMPA instances is the vast amount of samples required, which may be improved using more sophisticated neural network architectures or training algorithms. The DRL approach can be extended to also optimize $K$-DTMPA instances for other performance metrics besides the discounted cost criterion, e.g., the average cost criterion.

\small
\paragraph{Acknowledgements} We would like to thank Verus Pronk and Jan Korst for the discussions that formed the basis for this paper. This work used the Dutch national e-infrastructure with the support of the
SURF Cooperative using grant no. EINF-5192. This work was supported by the Netherlands Organisation for Scientific Research (NWO). Project: NWO Big data - Real Time ICT for Logistics. Number: 628.009.012. The work of Stella Kapodistria is supported by NWO through the Gravitation-grant NETWORKS-024.002.003. 
%load bibliography
\bibliographystyle{elsarticle-harv}
\bibliography{references.bib}

\begin{thebibliography}{34}
\expandafter\ifx\csname natexlab\endcsname\relax\def\natexlab#1{#1}\fi
\providecommand{\url}[1]{\texttt{#1}}
\providecommand{\href}[2]{#2}
\providecommand{\path}[1]{#1}
\providecommand{\DOIprefix}{doi:}
\providecommand{\ArXivprefix}{arXiv:}
\providecommand{\URLprefix}{URL: }
\providecommand{\Pubmedprefix}{pmid:}
\providecommand{\doi}[1]{\href{http://dx.doi.org/#1}{\path{#1}}}
\providecommand{\Pubmed}[1]{\href{pmid:#1}{\path{#1}}}
\providecommand{\bibinfo}[2]{#2}
\ifx\xfnm\relax \def\xfnm[#1]{\unskip,\space#1}\fi
%Type = Article
\bibitem[{Abdul-Malak and Kharoufeh(2018)}]{abdul2018optimally}
\bibinfo{author}{Abdul-Malak, D.T.}, \bibinfo{author}{Kharoufeh, J.P.}, \bibinfo{year}{2018}.
\newblock \bibinfo{title}{Optimally replacing multiple systems in a shared environment}.
\newblock \bibinfo{journal}{Probability in the Engineering and Informational Sciences} \bibinfo{volume}{32}, \bibinfo{pages}{179--206}.
%Type = Article
\bibitem[{Afrati et~al.(1986)Afrati, Cosmadakis, Papadimitriou, Papageorgiou and Papakostantinou}]{afrati1986complexity}
\bibinfo{author}{Afrati, F.}, \bibinfo{author}{Cosmadakis, S.}, \bibinfo{author}{Papadimitriou, C.H.}, \bibinfo{author}{Papageorgiou, G.}, \bibinfo{author}{Papakostantinou, N.}, \bibinfo{year}{1986}.
\newblock \bibinfo{title}{The complexity of the travelling repairman problem}.
\newblock \bibinfo{journal}{RAIRO - Theoretical Informatics and Applications - Informatique Th{\'e}orique et Applications} \bibinfo{volume}{20}, \bibinfo{pages}{79--87}.
%Type = Article
\bibitem[{Akcay(2022)}]{AKCAY2021}
\bibinfo{author}{Akcay, A.}, \bibinfo{year}{2022}.
\newblock \bibinfo{title}{An alert-assisted inspection policy for a production process with imperfect condition signals}.
\newblock \bibinfo{journal}{European Journal of Operational Research} \bibinfo{volume}{298}, \bibinfo{pages}{510--525}.
%Type = Techreport
\bibitem[{Bertsimas and Van~Ryzin(1989)}]{bertsimas1989dynamic}
\bibinfo{author}{Bertsimas, D.}, \bibinfo{author}{Van~Ryzin, G.}, \bibinfo{year}{1989}.
\newblock \bibinfo{title}{The dynamic traveling repairman problem}.
\newblock \bibinfo{type}{Working Paper} \bibinfo{number}{3036-89-MS}. MIT Sloan School of Management.
%Type = Article
\bibitem[{Boute et~al.(2022)Boute, Gijsbrechts, {van Jaarsveld} and Vanvuchelen}]{boute2021deep}
\bibinfo{author}{Boute, R.N.}, \bibinfo{author}{Gijsbrechts, J.}, \bibinfo{author}{{van Jaarsveld}, W.}, \bibinfo{author}{Vanvuchelen, N.}, \bibinfo{year}{2022}.
\newblock \bibinfo{title}{Deep reinforcement learning for inventory control: A roadmap}.
\newblock \bibinfo{journal}{European Journal of Operational Research} \bibinfo{volume}{298}, \bibinfo{pages}{401--412}.
%Type = Article
\bibitem[{Camci(2014)}]{camci2014travelling}
\bibinfo{author}{Camci, F.}, \bibinfo{year}{2014}.
\newblock \bibinfo{title}{The travelling maintainer problem: integration of condition-based maintenance with the travelling salesman problem}.
\newblock \bibinfo{journal}{Journal of the Operational Research Society} \bibinfo{volume}{65}, \bibinfo{pages}{1423--1436}.
%Type = Article
\bibitem[{Camci(2015)}]{camci2015maintenance}
\bibinfo{author}{Camci, F.}, \bibinfo{year}{2015}.
\newblock \bibinfo{title}{Maintenance scheduling of geographically distributed assets with prognostics information}.
\newblock \bibinfo{journal}{European Journal of Operational Research} \bibinfo{volume}{245}, \bibinfo{pages}{506--516}.
%Type = Article
\bibitem[{da~Costa et~al.(2023)da~Costa, Verleijsdonk, Voorberg, Akcay, Kapodistria, van Jaarsveld and Zhang}]{dtmpa}
\bibinfo{author}{da~Costa, P.}, \bibinfo{author}{Verleijsdonk, P.}, \bibinfo{author}{Voorberg, S.}, \bibinfo{author}{Akcay, A.}, \bibinfo{author}{Kapodistria, S.}, \bibinfo{author}{van Jaarsveld, W.}, \bibinfo{author}{Zhang, Y.}, \bibinfo{year}{2023}.
\newblock \bibinfo{title}{Policies for the dynamic traveling maintainer problem with alerts}.
\newblock \bibinfo{journal}{European Journal of Operational Research} \bibinfo{volume}{305}, \bibinfo{pages}{1141--1152}.
%Type = Article
\bibitem[{De~Jonge et~al.(2016)De~Jonge, Klingenberg, Teunter and Tinga}]{DEJONGE201693}
\bibinfo{author}{De~Jonge, B.}, \bibinfo{author}{Klingenberg, W.}, \bibinfo{author}{Teunter, R.}, \bibinfo{author}{Tinga, T.}, \bibinfo{year}{2016}.
\newblock \bibinfo{title}{Reducing costs by clustering maintenance activities for multiple critical units}.
\newblock \bibinfo{journal}{Reliability Engineering \& System Safety} \bibinfo{volume}{145}, \bibinfo{pages}{93--103}.
%Type = Article
\bibitem[{De~Jonge and Scarf(2020)}]{DEJONGE2020805}
\bibinfo{author}{De~Jonge, B.}, \bibinfo{author}{Scarf, P.A.}, \bibinfo{year}{2020}.
\newblock \bibinfo{title}{A review on maintenance optimization}.
\newblock \bibinfo{journal}{European Journal of Operational Research} \bibinfo{volume}{285}, \bibinfo{pages}{805--824}.
%Type = Article
\bibitem[{{De Moor} et~al.(2022){De Moor}, Gijsbrechts and Boute}]{DEMOOR2022535}
\bibinfo{author}{{De Moor}, B.J.}, \bibinfo{author}{Gijsbrechts, J.}, \bibinfo{author}{Boute, R.N.}, \bibinfo{year}{2022}.
\newblock \bibinfo{title}{Reward shaping to improve the performance of deep reinforcement learning in perishable inventory management}.
\newblock \bibinfo{journal}{European Journal of Operational Research} \bibinfo{volume}{301}, \bibinfo{pages}{535--545}.
%Type = Article
\bibitem[{Derman(1963)}]{derman1963optimal}
\bibinfo{author}{Derman, C.}, \bibinfo{year}{1963}.
\newblock \bibinfo{title}{On optimal replacement rules when changes of state are markovian}.
\newblock \bibinfo{journal}{Mathematical Optimization Techniques} \bibinfo{volume}{396}, \bibinfo{pages}{201--210}.
%Type = Article
\bibitem[{Drent et~al.(2020)Drent, Keizer and van Houtum}]{drent2020dynamic}
\bibinfo{author}{Drent, C.}, \bibinfo{author}{Keizer, M.O.}, \bibinfo{author}{van Houtum, G.J.}, \bibinfo{year}{2020}.
\newblock \bibinfo{title}{Dynamic dispatching and repositioning policies for fast-response service networks}.
\newblock \bibinfo{journal}{European Journal of Operational Research} \bibinfo{volume}{285}, \bibinfo{pages}{583--598}.
%Type = Inproceedings
\bibitem[{Feng et~al.(2021)Feng, Gluzman and Dai}]{ride_hailing}
\bibinfo{author}{Feng, J.}, \bibinfo{author}{Gluzman, M.}, \bibinfo{author}{Dai, J.G.}, \bibinfo{year}{2021}.
\newblock \bibinfo{title}{Scalable deep reinforcement learning for ride-hailing}, in: \bibinfo{booktitle}{2021 American Control Conference}, pp. \bibinfo{pages}{3743--3748}.
%Type = Article
\bibitem[{Gijsbrechts et~al.(2022)Gijsbrechts, Boute, Van~Mieghem and Zhang}]{inventory_management}
\bibinfo{author}{Gijsbrechts, J.}, \bibinfo{author}{Boute, R.N.}, \bibinfo{author}{Van~Mieghem, J.A.}, \bibinfo{author}{Zhang, D.J.}, \bibinfo{year}{2022}.
\newblock \bibinfo{title}{Can deep reinforcement learning improve inventory management? performance on lost sales, dual-sourcing, and multi-echelon problems}.
\newblock \bibinfo{journal}{Manufacturing \& Service Operations Management} \bibinfo{volume}{24}, \bibinfo{pages}{1349--1368}.
%Type = Article
\bibitem[{Haviv and Puterman(1992)}]{HAVIV1992267}
\bibinfo{author}{Haviv, M.}, \bibinfo{author}{Puterman, M.L.}, \bibinfo{year}{1992}.
\newblock \bibinfo{title}{Estimating the value of a discounted reward process}.
\newblock \bibinfo{journal}{Operations Research Letters} \bibinfo{volume}{11}, \bibinfo{pages}{267--272}.
%Type = Inproceedings
\bibitem[{Holler et~al.(2019)Holler, Vuorio, Qin, Tang, Jiao, Jin, Singh, Wang and Ye}]{holler2019deep}
\bibinfo{author}{Holler, J.}, \bibinfo{author}{Vuorio, R.}, \bibinfo{author}{Qin, Z.}, \bibinfo{author}{Tang, X.}, \bibinfo{author}{Jiao, Y.}, \bibinfo{author}{Jin, T.}, \bibinfo{author}{Singh, S.}, \bibinfo{author}{Wang, C.}, \bibinfo{author}{Ye, J.}, \bibinfo{year}{2019}.
\newblock \bibinfo{title}{Deep reinforcement learning for multi-driver vehicle dispatching and repositioning problem}, in: \bibinfo{booktitle}{2019 IEEE International Conference on Data Mining}, \bibinfo{organization}{IEEE}. pp. \bibinfo{pages}{1090--1095}.
%Type = Article
\bibitem[{Jagtenberg et~al.(2017)Jagtenberg, Bhulai and van~der Mei}]{jagtenberg}
\bibinfo{author}{Jagtenberg, C.J.}, \bibinfo{author}{Bhulai, S.}, \bibinfo{author}{van~der Mei, R.D.}, \bibinfo{year}{2017}.
\newblock \bibinfo{title}{Dynamic ambulance dispatching: is the closest-idle policy always optimal?}
\newblock \bibinfo{journal}{Health Care Management Science} \bibinfo{volume}{20}, \bibinfo{pages}{517--531}.
%Type = Article
\bibitem[{Ji et~al.(2019)Ji, Zheng, Wang and Li}]{ji2019deep}
\bibinfo{author}{Ji, S.}, \bibinfo{author}{Zheng, Y.}, \bibinfo{author}{Wang, Z.}, \bibinfo{author}{Li, T.}, \bibinfo{year}{2019}.
\newblock \bibinfo{title}{A deep reinforcement learning-enabled dynamic redeployment system for mobile ambulances}.
\newblock \bibinfo{journal}{Proceedings of the ACM on Interactive, Mobile, Wearable and Ubiquitous Technologies} \bibinfo{volume}{3}, \bibinfo{pages}{1--20}.
%Type = Article
\bibitem[{Keizer et~al.(2017)Keizer, Flapper and Teunter}]{OLDEKEIZER2017405}
\bibinfo{author}{Keizer, M.C.O.}, \bibinfo{author}{Flapper, S.D.P.}, \bibinfo{author}{Teunter, R.H.}, \bibinfo{year}{2017}.
\newblock \bibinfo{title}{Condition-based maintenance policies for systems with multiple dependent components: A review}.
\newblock \bibinfo{journal}{European Journal of Operational Research} \bibinfo{volume}{261}, \bibinfo{pages}{405--420}.
%Type = Misc
\bibitem[{Khorasgani et~al.(2020)Khorasgani, Wang and Gupta}]{khorasgani2020challenges}
\bibinfo{author}{Khorasgani, H.}, \bibinfo{author}{Wang, H.}, \bibinfo{author}{Gupta, C.}, \bibinfo{year}{2020}.
\newblock \bibinfo{title}{Challenges of applying deep reinforcement learning in dynamic dispatching}.
\newblock \href{http://arxiv.org/abs/2011.05570}{{\tt arXiv:2011.05570}}.
%Type = Misc
\bibitem[{Mnih et~al.(2013)Mnih, Kavukcuoglu, Silver, Graves, Antonoglou, Wierstra and Riedmiller}]{atari}
\bibinfo{author}{Mnih, V.}, \bibinfo{author}{Kavukcuoglu, K.}, \bibinfo{author}{Silver, D.}, \bibinfo{author}{Graves, A.}, \bibinfo{author}{Antonoglou, I.}, \bibinfo{author}{Wierstra, D.}, \bibinfo{author}{Riedmiller, M.}, \bibinfo{year}{2013}.
\newblock \bibinfo{title}{Playing atari with deep reinforcement learning}.
\newblock \href{http://arxiv.org/abs/1312.5602}{{\tt arXiv:1312.5602}}.
%Type = Misc
\bibitem[{Pechina et~al.(2019)Pechina, Usanov, van~de Ven and van~der Mei}]{pechina2019real}
\bibinfo{author}{Pechina, A.}, \bibinfo{author}{Usanov, D.}, \bibinfo{author}{van~de Ven, P.}, \bibinfo{author}{van~der Mei, R.}, \bibinfo{year}{2019}.
\newblock \bibinfo{title}{Real-time dispatching and relocation of emergency service engineers}.
\newblock \href{http://arxiv.org/abs/1910.01427}{{\tt arXiv:1910.01427}}.
%Type = Article
\bibitem[{Powell(2019)}]{powell2019unified}
\bibinfo{author}{Powell, W.B.}, \bibinfo{year}{2019}.
\newblock \bibinfo{title}{A unified framework for stochastic optimization}.
\newblock \bibinfo{journal}{European Journal of Operational Research} \bibinfo{volume}{275}, \bibinfo{pages}{795--821}.
%Type = Inproceedings
\bibitem[{Refaei~Afshar et~al.(2020)Refaei~Afshar, Zhang, Firat and Kaymak}]{StateAggregation}
\bibinfo{author}{Refaei~Afshar, R.}, \bibinfo{author}{Zhang, Y.}, \bibinfo{author}{Firat, M.}, \bibinfo{author}{Kaymak, U.}, \bibinfo{year}{2020}.
\newblock \bibinfo{title}{A state aggregation approach for solving knapsack problem with deep reinforcement learning}, in: \bibinfo{editor}{Pan, S.J.}, \bibinfo{editor}{Sugiyama, M.} (Eds.), \bibinfo{booktitle}{Proceedings of The 12th Asian Conference on Machine Learning}, \bibinfo{publisher}{PMLR}. pp. \bibinfo{pages}{81--96}.
%Type = Article
\bibitem[{Schmid(2012)}]{schmid2012solving}
\bibinfo{author}{Schmid, V.}, \bibinfo{year}{2012}.
\newblock \bibinfo{title}{Solving the dynamic ambulance relocation and dispatching problem using approximate dynamic programming}.
\newblock \bibinfo{journal}{European Journal of Operational Research} \bibinfo{volume}{219}, \bibinfo{pages}{611--621}.
%Type = Book
\bibitem[{Schrijver(2003)}]{schrijver2003combinatorial}
\bibinfo{author}{Schrijver, A.}, \bibinfo{year}{2003}.
\newblock \bibinfo{title}{Combinatorial optimization: polyhedra and efficiency}. volume~\bibinfo{volume}{24}.
\newblock \bibinfo{edition}{2} ed., \bibinfo{publisher}{Springer}.
%Type = Article
\bibitem[{Silver et~al.(2018)Silver, Hubert, Schrittwieser, Antonoglou, Lai, Guez, Lanctot, Sifre, Kumaran, Graepel, Lillicrap, Simonyan and Hassabis}]{chess}
\bibinfo{author}{Silver, D.}, \bibinfo{author}{Hubert, T.}, \bibinfo{author}{Schrittwieser, J.}, \bibinfo{author}{Antonoglou, I.}, \bibinfo{author}{Lai, M.}, \bibinfo{author}{Guez, A.}, \bibinfo{author}{Lanctot, M.}, \bibinfo{author}{Sifre, L.}, \bibinfo{author}{Kumaran, D.}, \bibinfo{author}{Graepel, T.}, \bibinfo{author}{Lillicrap, T.}, \bibinfo{author}{Simonyan, K.}, \bibinfo{author}{Hassabis, D.}, \bibinfo{year}{2018}.
\newblock \bibinfo{title}{A general reinforcement learning algorithm that masters chess, shogi, and go through self-play}.
\newblock \bibinfo{journal}{Science} \bibinfo{volume}{362}, \bibinfo{pages}{1140--1144}.
%Type = Misc
\bibitem[{Snellius(accessed: 09.18.2022)}]{snellius}
\bibinfo{author}{Snellius}, \bibinfo{year}{accessed: 09.18.2022}.
\newblock \bibinfo{title}{Snellius {S}upercomputer}.
\newblock \URLprefix \url{https://www.surf.nl/en/dutch-national-supercomputer-snellius}.
%Type = Misc
\bibitem[{{SURF}(accessed: 12.06.2023)}]{surf2022}
\bibinfo{author}{{SURF}}, \bibinfo{year}{accessed: 12.06.2023}.
\newblock \bibinfo{title}{Surf services and rates 2023}.
\newblock \URLprefix \url{https://www.surf.nl/files/2023-08/surf-services-and-rates-2023_version-aug-2022.pdf}.
%Type = Misc
\bibitem[{Temizöz et~al.(2023)Temizöz, Imdahl, Dijkman, Lamghari-Idrissi and van Jaarsveld}]{temizöz2023deep}
\bibinfo{author}{Temizöz, T.}, \bibinfo{author}{Imdahl, C.}, \bibinfo{author}{Dijkman, R.}, \bibinfo{author}{Lamghari-Idrissi, D.}, \bibinfo{author}{van Jaarsveld, W.}, \bibinfo{year}{2023}.
\newblock \bibinfo{title}{Deep controlled learning for inventory control}.
\newblock \href{http://arxiv.org/abs/2011.15122v8}{{\tt arXiv:2011.15122v8}}.
%Type = Article
\bibitem[{Van~Buuren et~al.(2018)Van~Buuren, Jagtenberg, Van~Barneveld, Van Der~Mei and Bhulai}]{jagtenberg2}
\bibinfo{author}{Van~Buuren, M.}, \bibinfo{author}{Jagtenberg, C.}, \bibinfo{author}{Van~Barneveld, T.}, \bibinfo{author}{Van Der~Mei, R.}, \bibinfo{author}{Bhulai, S.}, \bibinfo{year}{2018}.
\newblock \bibinfo{title}{Ambulance dispatch center pilots proactive relocation policies to enhance effectiveness}.
\newblock \bibinfo{journal}{Interfaces} \bibinfo{volume}{48}, \bibinfo{pages}{235--246}.
%Type = Unpublished
\bibitem[{Vanvuchelen et~al.(2022)Vanvuchelen, De~Moor and Boute}]{Vanvuchelen2022}
\bibinfo{author}{Vanvuchelen, N.}, \bibinfo{author}{De~Moor, B.J.}, \bibinfo{author}{Boute, R.N.}, \bibinfo{year}{2022}.
\newblock \bibinfo{title}{The use of continuous action representations to scale deep reinforcement learning: An application to inventory control}.
\newblock \bibinfo{note}{Available at SSRN: \url{https://ssrn.com/abstract=4253600}}.
%Type = Inproceedings
\bibitem[{Zhang et~al.(2020)Zhang, Odonkor, Zheng, Khorasgani, Serita, Gupta and Wang}]{zhang2020dynamic}
\bibinfo{author}{Zhang, C.}, \bibinfo{author}{Odonkor, P.}, \bibinfo{author}{Zheng, S.}, \bibinfo{author}{Khorasgani, H.}, \bibinfo{author}{Serita, S.}, \bibinfo{author}{Gupta, C.}, \bibinfo{author}{Wang, H.}, \bibinfo{year}{2020}.
\newblock \bibinfo{title}{Dynamic dispatching for large-scale heterogeneous fleet via multi-agent deep reinforcement learning}, in: \bibinfo{booktitle}{2020 IEEE International Conference on Big Data}, \bibinfo{organization}{IEEE}. pp. \bibinfo{pages}{1436--1441}.

\end{thebibliography}
\pagebreak
%load appendix

\section*{Supplementary material}
\beginsupplement
\appendix

\section{Degradation matrices}\label{app: matrices}
The following degradation matrices are adopted from \citet[Section 6.1]{dtmpa}.
\[\small{
\textrm{Q2} = \begin{bmatrix}
0.8 & 0.2 & 0 & 0 & 0\\
0 & 0.7 & 0.3 & 0 & 0\\
0 & 0 & 0.7 & 0.3 & 0 \\
0 & 0 & 0 & 0.7 & 0.3 \\
0 & 0 & 0 & 0 & 1 
\end{bmatrix}
}
\]
\\
\[\small{
\textrm{Q3} = \begin{bmatrix}
0.8 & 0.2 & 0 & 0 & 0\\
0 & 0.3 & 0.7 & 0 & 0\\
0 & 0 & 0.3 & 0.7 & 0 \\
0 & 0 & 0 & 0.3 & 0.7 \\
0 & 0 & 0 & 0 & 1 
\end{bmatrix}\quad
\textrm{Q4} = \begin{bmatrix}
0.8 & 0.2 & 0 & 0 & 0 & 0 & 0\\
0 & 0.7 & 0.3 & 0 & 0 & 0 & 0 \\
0 & 0 & 0.7 & 0.3 & 0 & 0 & 0 \\
0 & 0 & 0 & 0.7 & 0.3 & 0 & 0 \\
0 & 0 & 0 & 0 & 0.7 & 0.3 & 0 \\
0 & 0 & 0 & 0 & 0 & 0.7 & 0.3 \\
0 & 0 & 0 & 0 & 0 & 0 & 1 \\
\end{bmatrix}}
\]
\\

\newpage
\section{Nomenclature} \label{app:nomenc}
\begin{table*}[ht]
 \begin{framed}
 \printnomenclature
 \nomenclature[01]{$t$}{Time step}
 \nomenclature[02]{$\mathcal{M}$}{Set of machines}
 \nomenclature[03]{$\mathcal{K}$}{Set of engineers}
 \nomenclature[04]{$M$}{Number of machines}
 \nomenclature[05]{$m$}{Machine index}
 \nomenclature[06]{$K$}{Number of engineers}
 \nomenclature[07]{$k$}{Engineer index}
 \nomenclature[08]{$c^\textrm{PM}_m$}{PM cost of machine $m$}
 \nomenclature[09]{$c^\textrm{CM}_m$}{CM cost of machine $m$} 
 \nomenclature[10]{$c^\textrm{DT}_m$}{Downtime cost of machine $m$} 
 \nomenclature[11]{$c^\textrm{T}$}{Travel cost} 
\nomenclature[12]{$\gamma$}{Discount factor} 
 \nomenclature[13]{$\theta_{ij}$}{Travel time between machines $i$ and $j$}
 \nomenclature[14]{$x^\textbf{h}_m$}{As-good-as-new state}
 \nomenclature[15]{$x^\textbf{\textbf{f}}_m$}{Failed state}
 \nomenclature[16]{$\mathcal{N}_m$}{State space of machine $m$}
 \nomenclature[17]{$x_m(t)$}{Degradation process of machine $m$} 
 \nomenclature[18]{$T^{x_m}_m$}{Random transition time from state $x_m$ to $x_m+1$}
\nomenclature[19]{$h$}{Network state}
\nomenclature[20]{$f(h)$}{Feature vector}
\nomenclature[21]{$\mathcal{U}(h)$}{State-dependent action set}
\nomenclature[22]{$C(h,a)$}{Cost function}
\nomenclature[23]{$\pi$}{Policy}
 % ADD THE FOLLOWING LINE IF THE NUMBER OF ITEMS IS ODD
 %\nomenclature{\mbox{}}{}
 \end{framed}
\end{table*}

\newpage
\section{Deep reinforcement learning hyperparameters}\label{app: drl_param}
We list the choice of hyperparameters for the training algorithm per $K$-DTMPA instance. All neural networks utilize the rectified linear unit (ReLU) activation function.

\begin{table}[!ht]
\centering
\begin{tabular}{c|l|c}
\toprule
Hyperparameter & Description & Value \\ \midrule
$L$ & number of neural network layers & $3$ \\
$d^1$ & dimension of input layer & $128$ \\
$d^l$ & dimension of hidden layer $l=2, \ldots, L$ & $64$ \\
$\text{MAX\_SAMPLES}$ & number of samples & 150,000 \\
$r_\textrm{min}$ & minimum number of roll-outs & $1500$ \\
$r_\textrm{max}$ & maximum number of roll-outs & $7500$ \\
$\text{BATCH\_SIZE}$ & batch size & $64$\\
$k$ & $k$-value bandit optimizer & $2.0$ \\
$\epsilon$ & fraction random actions & $0.02$ \\
 \bottomrule
\end{tabular}%
\caption{Approximate policy iteration hyperparameters for all single maintainer instances, i.e., M4K1-Q2Q3C2, M6K1-Q2Q3Q4C2, and the $1$-DTMPA instances induced by the clusters when training the decomposition heuristic.}
\label{tab: hyperM4K1-Q2Q3}
\end{table}

\begin{table}[!ht]
\centering
\begin{tabular}{c|l|c}
\toprule
Hyperparameter & Description & Value \\ \midrule
$L$ & number of neural network layers & $4$ \\
$d^1$ & dimension of input layer & $256$ \\
$d^l$ & dimension of hidden layer $l=2, \ldots, L$ & $128$ \\
$\text{MAX\_SAMPLES}$ & number of samples & 500,000 \\
$r_\textrm{min}$ & minimum number of roll-outs & $1500$ \\
$r_\textrm{max}$ & maximum number of roll-outs & $7500$ 
\\
$\text{BATCH\_SIZE}$ & batch size & $64$\\
$k$ & $k$-value bandit optimizer & $2.0$ \\
$\epsilon$ & fraction random actions & $0.02$ \\
 \bottomrule
\end{tabular}%
\caption{Approximate policy iteration hyperparameters for the dispatching \& repositioning instance M8K3-$\tilde{\textrm{Q}}$1C1.}
\label{tab: hyperM8K3-Q1C1}
\end{table}

\begin{table}[!ht]
\centering
\begin{tabular}{c|l|c}
\toprule
Hyperparameter & Description & Value \\ \midrule
$L$ & number of neural network layers & $4$ \\
$d^1$ & dimension of input layer & $256$ \\
$d^l$ & dimension of hidden layer $l=2, \ldots, L$ & $128$ \\
$\text{MAX\_SAMPLES}$ & number of samples & 850,000 \\
$r_\textrm{min}$ & minimum number of roll-outs & $1500$ \\
$r_\textrm{max}$ & maximum number of roll-outs & $7500$ 
\\
$\text{BATCH\_SIZE}$ & batch size & $64$\\
$k$ & $k$-value bandit optimizer & $2.0$ \\
$\epsilon$ & fraction random actions & $0.02$ \\
 \bottomrule
\end{tabular}
\caption{Approximate policy iteration hyperparameters for the preventive maintenance instance M8K3-$\tilde{\textrm{Q}}$2C3.}
\label{tab: hyperM8K3-Q2C3}
\end{table}

\begin{table}[!ht]
\centering
\begin{tabular}{c|l|c}
\toprule
Hyperparameter & Description & Value \\ \midrule
$L$ & number of neural network layers & $4$ \\
$d^1$ & dimension of input layer & $512$ \\
$d^l$ & dimension of hidden layer $l=2, \ldots, L$ & $256$ \\
$\text{MAX\_SAMPLES}$ & number of samples & 2,000,000 \\
$r_\textrm{min}$ & minimum number of roll-outs & $(\text{gen1:}1500, \text{gen2:}500)$ \\
$r_\textrm{max}$ & maximum number of roll-outs & $(\text{gen1:}7500, \text{gen2:}1500)$ 
\\
$\text{BATCH\_SIZE}$ & batch size & $64$\\
$k$ & $k$-value bandit optimizer & $2.0$ \\
$\epsilon$ & fraction random actions & $0.02$ \\
 \bottomrule
\end{tabular}
\caption{Approximate policy iteration hyperparameters for the dispatching \& repositioning instance M35K5-$\tilde{\textrm{Q}}$3C1. Note that the data for the second generation is collected using less roll-outs.}
\label{tab: hyperM35K5-Q3C1_step1}
\end{table}

\begin{table}[!ht]
\centering
\begin{tabular}{c|l|c}
\toprule
Hyperparameter & Description & Value \\ \midrule
$L$ & number of neural network layers & $4$ \\
$d^1$ & dimension of input layer & $512$ \\
$d^l$ & dimension of hidden layer $l=2, \ldots, L$ & $256$ \\
$\text{MAX\_SAMPLES}$ & number of samples & 3,000,000 \\
$r_\textrm{min}$ & minimum number of roll-outs & $500$ \\
$r_\textrm{max}$ & maximum number of roll-outs & $1500$ 
\\
$\text{BATCH\_SIZE}$ & batch size & $64$\\
$k$ & $k$-value bandit optimizer & $2.0$ \\
$\epsilon$ & fraction random actions & $0.02$ \\
 \bottomrule
\end{tabular}
\caption{Approximate policy iteration hyperparameters for the preventive maintenance instance M35K5-$\tilde{\textrm{Q}}$4C3.}
\label{tab: hyperM35K5-Q3C1_step2}
\end{table}

\clearpage
\section{Deep reinforcement learning training durations and costs}\label{app: drl_training}
We list some information regarding the cost and duration of the sample collection for the $K$-DTMPA instances M8K3-$\tilde{\textrm{Q}}$1C1, M8K3-$\tilde{\textrm{Q}}$2C3, M35K5-$\tilde{\textrm{Q}}$3C1 and M35K5-$\tilde{\textrm{Q}}$4C3.

\begin{table}[H]
\centering
\begin{tabular}{c|c|c|c|c}
\toprule
 & M8K3-$\tilde{\textrm{Q}}$1C1 & M8K3-$\tilde{\textrm{Q}}$2C3 & M35K5-$\tilde{\textrm{Q}}$3C1 & M35K5-$\tilde{\textrm{Q}}$4C3 \\ \midrule
Training time (gen 1) & 0.59 hrs & 0.88 hrs & 13.67 hrs & 6.45 hrs\\ 
%Cost data gen 1: SBU & 755 SBU & 1126 SBU & 17497 SBU & 8256 SBU \\ 
Cost & \euro7.55 & \euro11.26 & \euro174.97 & \euro82.56 \\ \midrule
Training time (gen 2) & 2.17 hrs & 3.39 hrs & 65.23 hrs & - \\ 
%Cost data gen 2: SBU & 2777 SBU & 4339 SBU & 83494 SBU & - \\ 
Cost & \euro27.77 & \euro43.39 & \euro834.94 & -
\end{tabular}
\caption{Cost and duration of the sample collection phase. Each reported time is a duration estimate to collect the required amount of samples (cf. \ref{app: drl_param}) using 10 \emph{thin computing nodes} on the Dutch National Supercomputer \cite{snellius}. The corresponding cost is computed from \citet{surf2022}.}
\label{tab: trainingtimes}
\end{table}
%% \bibitem{key}...
% \end{thebibliography}
\end{document}